\documentclass[12pt]{article}

\usepackage[T1]{fontenc}

\usepackage{amsfonts}
\usepackage{latexsym}
\usepackage{color}
\usepackage{graphicx}
\usepackage{algorithmic}

\newtheorem{thm}{Theorem}[section]
\newtheorem{cor}{Corollary}[section]
\newtheorem{lem}{Lemma}[section]

\newtheorem{prob}{Problem}[section]
\newtheorem{exa}{Example}[section]

\newtheorem{defn}{Definition}[section]

\newcommand{\bfc}{\mbox{$\mbox{\boldmath $c$}$}} 
\newcommand{\bfe}{\mbox{$\mbox{\boldmath $e$}$}}
\newcommand{\bff}{\mbox{$\mbox{\boldmath $f$}$}}
\newcommand{\bfg}{\mbox{$\mbox{\boldmath $g$}$}} 
\newcommand{\bfh}{\mbox{$\mbox{\boldmath $h$}$}} 
\newcommand{\bfk}{\mbox{$\mbox{\boldmath $k$}$}} 
\newcommand{\bfp}{\mbox{$\mbox{\boldmath $p$}$}} 
\newcommand{\bfq}{\mbox{$\mbox{\boldmath $q$}$}} 
\newcommand{\bfr}{\mbox{$\mbox{\boldmath $r$}$}} 
\newcommand{\bfu}{\mbox{$\mbox{\boldmath $u$}$}} 
\newcommand{\bfv}{\mbox{$\mbox{\boldmath $v$}$}} 
\newcommand{\bfx}{\mbox{$\mbox{\boldmath $x$}$}} 
\newcommand{\bfy}{\mbox{$\mbox{\boldmath $y$}$}} 

\newcommand{\sbff}{\mbox{\boldmath \scriptsize $f$}} 
\newcommand{\sbfg}{\mbox{\boldmath \scriptsize $g$}} 
\newcommand{\sbfh}{\mbox{\boldmath \scriptsize $h$}} 
\newcommand{\sbfk}{\mbox{\boldmath \scriptsize $k$}} 
\newcommand{\sbfp}{\mbox{\boldmath \scriptsize $p$}} 
\newcommand{\sbfq}{\mbox{\boldmath \scriptsize $q$}} 
\newcommand{\sbfr}{\mbox{\boldmath \scriptsize $r$}} 
\newcommand{\sbfv}{\mbox{\boldmath \scriptsize $v$}} 

\newcommand{\bphi}{\mbox{$\mbox{\boldmath $\varphi$}$}} 
\newcommand{\bpsi}{\mbox{$\mbox{\boldmath $\psi$}$}} 
\newcommand{\bmu}{\mbox{$\mbox{\boldmath $\mu$}$}} 
\newcommand{\bomega}{\mbox{$\mbox{\boldmath $\omega$}$}} 

\newcommand{\bfzero}{\mbox{\boldmath $0$}}

\newcommand{\sbfomega}{\mbox{\boldmath \scriptsize $\omega$}} 

\newcommand{\dlangle}{\mbox{$\mbox{$\langle\! \langle$}$}}
\newcommand{\drangle}{\mbox{$\mbox{$\rangle\! \rangle$}$}}
\newcommand{\dnorm}{\mbox{$\mbox{$|\hspace{-0.5mm}|\hspace{-0.5mm}|$}$}}

\begin{document}

\title{An optimal linear filter for estimation of random functions in Hilbert space}

\author{Phil Howlett and Anatoli Torokhti}

\maketitle

\begin{abstract}
Let $\bff$ be a square-integrable, zero-mean, random vector with observable realizations in a Hilbert space $H$, and let $\bfg$ be an associated square-integrable, zero-mean, random vector with realizations, which are not observable, in a Hilbert space $K$.  We seek an optimal filter in the form of a closed linear operator $X$ acting on the observable realizations of a proximate vector $\bff_{\epsilon} \approx \bff$ that provides the best estimate $\widehat{\bfg}_{\epsilon} = X\! \bff_{\epsilon}$ of the vector $\bfg$.  We assume the required covariance operators are known.  The results are illustrated with a typical example.
\end{abstract}



\section{Introduction}
\label{intro}

A common problem in engineering, applied mathematics and statistics is the estimation of a random function $\bfg$, whose realizations $\bfg(\bomega)$ are not observable, by using the observable realizations $\bff(\bomega)$ of an associated random function $\bff$.  We consider the following problem.

\begin{prob}
\label{p1}
Let $(\Omega, \Sigma, \mu)$ be a probability space, $H$ and $K$ Hilbert spaces, and $\bff \in L^2(\Omega,H)$ and $\bfg \in L^2(\Omega,K)$ square-integrable, zero-mean, random functions with respective observable and unobservable realizations $\bff(\bomega)$ and $\bfg(\bomega)$ for each outcome $\bomega \in \Omega$.  Find a closed, densely defined, linear operator $X:D(X) \subseteq H \rightarrow K$, a proximate observable function $\bff_{\epsilon}$ for each $\epsilon > 0$ with ${\mathbb E}[ \| \bff_{\epsilon} - \bff \|^2] < \epsilon$ and $\bff_{\epsilon}(\bomega) \in D(X)$ for $\mu$-almost all $\bomega \in \Omega$, and a corresponding estimate $\widehat{\bfg}_{\epsilon} = X\!\bff_{\epsilon}$ of the unobservable function, such that
\begin{equation}
\label{generror}
{\mathbb E}[ \| X\!\bff_{\epsilon} - \bfg \|^2] = \mbox{$\int_{\Omega}$} \| X\!\bff_{\epsilon}(\bomega) - \bfg(\bomega) \|^2 \mu(d \bomega)
\end{equation}
is minimized.  $\hfill \Box$
\end{prob}

For each outcome $\bomega \in \Omega$ the realization $\bfr_{\epsilon}(\bomega) = X\!\bff_{\epsilon}(\bomega) - \bfg(\bomega)$ of the error function is an element of the Hilbert space $K$.  The value $\| \bfr_{\epsilon}(\bomega) \|^2$ is the square of the magnitude of the pointwise error.  The estimated overall error ${\mathbb E}[ \| X\!\bff_{\epsilon} - \bfg \|^2]$ in (\ref{generror}) is the mean or expected value of the square of the magnitude of the pointwise error.  The proximate observable function $\bff_{\epsilon}$ must be close to the observable function $\bff$ in the sense that the mean square observation error ${\mathbb E}[ \| \bff_{\epsilon} - \bff \|^2]$ must be small.  The outcomes $\bff_{\epsilon}(\bomega)$ must lie in the domain space $D(X)$ of the operator $X$ for $\mu$-almost all $\bomega \in \Omega$.   The pointwise estimate $\widehat{\bfg}_{\epsilon} = X\!\bff_{\epsilon} \in L^2(\Omega,K)$ of the unobservable function is defined by $\widehat{\bfg}_{\epsilon}(\bomega) = X\!\bff_{\epsilon}(\bomega) = X \cdot \bff_{\epsilon}(\bomega)$ for each $\bomega \in \Omega$.  The linear operator $X$ does not depend on the approximation parameter $\epsilon$.  We use the terms {\em random vector} and {\em random function} interchangeably, but for the most part, prefer the latter.

We assume that the key covariance operators, the auto-covariance $E_{\sbff \! \sbff}$ and the cross-covariance $E_{\sbfg \sbff}$, are known bounded linear operators.  We expect to find a solution in the form $\widehat{\bfg} = X \! \bff$ where $X = E_{\sbfg \sbff} {E_{\sbff \! \sbff}}^{\dag}$ and ${E_{\sbff \! \sbff}}^{\dag}$ is the generalized inverse auto-covariance operator and this perception enables us to identify some critical issues.  The auto-covariance operator $E_{\sbff \! \sbff}$ is positive semi-definite, self-adjoint and compact.  Therefore the spectral set is reduced to a countable collection of real non-negative eigenvalues.  When there are an infinite number of positive eigenvalues the auto-covariance is not bounded below and the range space is not a closed subspace.  Therefore the generalized inverse auto-covariance ${E_{\sbff \! \sbff}}^{\dag}$ is an unbounded linear operator.  Consequently the proposed solution $X = E_{\sbfg \sbff} {E_{\sbff \! \sbff}}^{\dag}$ is also unbounded.  Now there are two specific issues that must be resolved.

In the first instance the usual justification for the solution assumes that the auto-covariance of the transformed observable function $X\! \bff$ is given by the formula $E_{X\! \sbff, X\! \sbff} = XE_{\sbff \! \sbff} X^*$.  The usual justification is no longer valid if the operator $X$ is unbounded.  The matter is resolved by writing $X = T{E_{\sbff \! \sbff}}^{\dag}$ where $T$ is a bounded linear operator and then using an alternative argument to find an optimal value for $T$.

In the second instance a solution in the form $\widehat{\bfg} = E_{\sbfg \sbff} {E_{\sbff \! \sbff}}^{\dag} \bff$ would require the observable function $\bff$ to lie in the domain of the unbounded operator ${E_{\sbff \! \sbff}}^{\dag}$.  This cannot be guaranteed.  The difficulty can be resolved by introducing a proximate observable function $\bff_{\epsilon}$, for each $\epsilon > 0$, which must lie in the domain of ${E_{\sbff \! \sbff}}^{\dag}$ but needs to be close to $\bff$ in the sense that the mean square error in the observed values satisfies ${\mathbb E}[ \| \bff_{\epsilon} - \bff \|^2] < \epsilon$.  The proposed solution now takes the form $\widehat{\bfg}_{\epsilon} = E_{\sbfg \sbff} {E_{\sbff \! \sbff}}^{\dag} \bff_{\epsilon}$.  This raises a further question.  How can we ensure that the operator\textemdash which does not depend on the approximation parameter\textemdash is still optimal for the proximate function?  The answer is found by taking the proximate function as a partial sum of the Fourier series for the observable function. 

\subsection{A basic formulation of the problem.}
\label{bfp}

Let $(\Omega ,\Sigma ,\mu )$ denote a probability space where $\Omega $ is the set of outcomes, $\Sigma $ a complete $\sigma$-field of measurable subsets $E\subseteq \Omega $ and $\mu : \Sigma \rightarrow [0,1]$ an associated probability measure on $\Sigma $, with $\mu (\Omega )=1$.   Each element $\bomega \in \Omega $ represents the outcome of an observation or experiment and each $E \in \Sigma $ is a set of outcomes, called an event. We say that the event $E$ has occurred if $\bomega \in E$.  Let $\bff = \bff(\bomega) \in {\mathbb C}^m$ and $\bfg = \bfg(\bomega) \in {\mathbb C}^n$ be complex-valued random vectors with zero mean.  That is, we assume
$$
\bmu_{\sbff} = \mbox{$\int_{\Omega}$} \bff(\bomega) \mu(d \bomega) = \bfzero \quad \mbox{and} \quad \bmu_{\sbfg} = \mbox{$\int_{\Omega}$}\, \bfg(\bomega) \mu(d \bomega) = \bfzero.
$$
We would like to estimate $\bff$ from a knowledge of $\bfg$.  One might postulate a linear relationship in the form $X \bff = \bfg + \bfr$ where $X \in {\mathbb C}^{n \times m}$ is an unknown matrix and $\bfr = \bfr(\bomega) \in {\mathbb C}^n$ is a random error vector which is independent of $\bff$ and has zero mean.   If so then after $k$ realizations one would obtain a system of equations
\begin{equation}
\label{le1}
X [ \bff_1,\ldots,\bff_k] = [\bfg_1,\ldots,\bfg_k] + [\bfr_1, \ldots, \bfr_k] \iff XF = G + R
\end{equation}
where we have written $\bff_j = \bff(\bomega_j)$,  $\bfg_j = \bfg(\bomega_j)$, $\bfr_j = \bfr(\bomega_j)$ where $\bomega_j$ is the outcome of the $j^{\mbox{\scriptsize th}}$ experiment and where $F = [\bff_1,\ldots,\bff_k] \in {\mathbb C}^{m \times k}$, $G = [\bfg_1,\ldots,\bfg_k] \in {\mathbb C}^{n \times k}$ and $R = [\bfr_1,\ldots,\bfr_k] \in {\mathbb C}^{n \times k}$.  In general, because we have merely postulated a linear relationship, one would not expect this equation to be satisfied exactly.  Thus we seek to minimize the mean-square error $\sum_{j=1}^k \| X\! \bff_j - \bfg_j \|^2 = \sum_{j=1}^k \|\bfr_j\|^2$.   Hence we solve the system
\begin{equation}
\label{le2}
XFF^* = GF^* \quad \Leftrightarrow \quad X \mbox{$\sum_{i=1}^k$} \bff_i {\bff_i}^* = \mbox{$\sum_{i=1}^k$} \bfg_i {\bff_i}^*.
\end{equation}
We make a probabilistic interpretation of this equation by noting that
\begin{equation}
\label{acm}
E_{\sbff \! \sbff} = {\mathbb E}[ \bff \bff^*] \sim \frac{1}{k} \mbox{$\sum_{i=1}^k$} \bff_i {\bff_i}^* = \frac{1}{k}\, FF^*
\end{equation}
and
\begin{equation}
\label{cm}
E_{\sbfg \sbff} = {\mathbb E}[ \bfg \bff^*] \sim  \frac{1}{k} \mbox{$\sum_{i=1}^k$} \bfg_i {\bff_i}^* = \frac{1}{k}\, GF^*
\end{equation}
where ${\mathbb E}$ is the expectation operator and $E_{\sbff \! \sbff} \in {\mathbb C}^{m \times m}$ and $E_{\sbfg \sbff} \in {\mathbb C}^{n \times m}$ are the standard auto-covariance and cross-covariance matrices for zero-mean vectors.  Thus we rewrite the equation for the best estimate of $X$ in the form
\begin{equation}
\label{le3}
X E_{\sbff \! \sbff} = E_{\sbfg \sbff}.
\end{equation}

\subsection{The definitive properties of the covariance matrices.}
\label{dpcm}

To extend the above analysis to random vectors in Hilbert space we must be able to define appropriate covariance operators.  Notice that 
$$
\langle E_{\sbff \! \sbff}\, \bfx, \bfu \rangle \sim \frac{1}{k} \langle FF^* \bfx, \bfu \rangle = \frac{1}{k}\, \mbox{$\sum_{j=1}^k$} \langle \bfu, \bff_j \rangle \langle \bff_j, \bfx \rangle \sim {\mathbb E}[ \langle \bfu, \bff \rangle \langle \bff, \bfx \rangle]
$$
for each $\bfx, \bfu \in {\mathbb C}^m$ and
$$
\langle E_{\sbfg \sbff}\, \bfx, \bfy \rangle \sim \frac{1}{k} \langle GF^*\bfx, \bfy \rangle = \frac{1}{k}\, \mbox{$\sum_{j=1}^k$} \langle \bfy, \bfg_j \rangle \langle \bff_j, \bfx \rangle \sim {\mathbb E}[ \langle \bfy, \bfg \rangle \langle \bff, \bfx \rangle ] 
$$
for each $\bfx \in {\mathbb C}^m$ and $\bfy \in {\mathbb C}^n$ and also that
$$
\mbox{tr}(E_{\sbff \! \sbff}) \sim \frac{1}{k}\, \mbox{tr}(FF^*) = \frac{1}{k} \mbox{tr}(F^*F) = \frac{1}{k} \mbox{$\sum_{i=1}^k$} \| \bff_i \|^2 \sim {\mathbb E}[ \| \bff \|^2 ]
$$
where ${\mathbb E}$ is the expectation operator.  By taking the limit as the number of independent realizations tends to infinity we obtain the basic theoretical relationships
\begin{eqnarray}
\langle E_{\sbff \! \sbff}\, \bfx, \bfu \rangle & = & {\mathbb E}[ \langle \bfu, \bff \rangle \langle \bff, \bfx \rangle], \label{covmatprop1} \\
\langle E_{\sbfg \sbff}\, \bfx, \bfy \rangle & = & {\mathbb E}[ \langle \bfy, \bfg \rangle \langle \bff, \bfx \rangle ], \label{covmatprop2} \\
\mbox{tr}(E_{\sbff \! \sbff}) & = & {\mathbb E}[ \| \bff \|^2] \label{covmatprop3}
\end{eqnarray}
for all $\bfx, \bfu \in {\mathbb C}^m$ and $\bfy \in {\mathbb C}^n$.  We will take these as the definitive properties of the covariance operators for the random vectors $\bff$ and $\bfg$.

\subsection{A typical application\textemdash input retrieval in a linear system.}
\label{taerf}

We illustrate our theoretical results by considering the problem of input retrieval in an infinite-dimensional linear system.  Our formal task is to find an optimal estimate of the system input from observations of the system output.  The input is a random function $\bfg$ which is represented as a Fourier series with random coefficients.  The output $\bff$ is a random function where each realization $\bff(\omega)$ of the output is uniquely determined by the corresponding realization $\bfg(\omega)$ of the input for some $\omega \in \Omega$.  We assume there is no independently generated noise to disrupt our observations of the output.  This makes no substantial difference to the methodology.  The introduction of noise simply decreases the accuracy of the estimation.  In our hypothetical example we consider a known system so that the required covariance operators $E_{\sbff\! \sbff}$ and $E_{\sbfg \sbff}$ are also known.  In practice it may be necessary to estimate these operators {\em a priori} in a controlled experiment.  Each observed output $\bff(\omega)$ is approximated by a truncated Fourier series $\bff_{\epsilon}(\omega) = \bff_m(\omega)$ for some fixed $m \in {\mathbb N}$ and the input is then estimated using the formula $\widehat{\bfg}_{\epsilon} = E_{\sbfg \sbff} {E_{\sbff\! \sbff}}^{\dag} \bff_{\epsilon} \Leftrightarrow \widehat{\bfg}_m = E_{\sbfg \sbff} {E_{\sbff\! \sbff}}^{\dag}  \bff_m$.  See also \cite{how1} for an application to input retrieval in finite-dimensional linear control systems and \cite[Section 8.4.1, pp 261\textendash 262]{avr1} for the extension of these ideas to infinite-dimensional systems.   

Our hypothetical example is a special case of a more general collection of {\em so-called} inverse problems.  See Cotter et al.~\cite{cot1} for an extended discussion of the underlying statistical theory of optimal estimation and a collection of particular inverse problems arising from data assimilation in fluid mechanics.  In each application  one assumes that the system evolves in a predominantly deterministic manner from some unknown initial configuration and that the evolution is monitored either directly or indirectly by observation of various output signals that may or may not be disrupted by random noise.  The objective is to make inference about the underlying velocity field. For problems without model error the inference is on the initial conditions.  For problems with model error the inference is on the initial conditions and on the driving noise process or, equivalently, on the entire time-dependent velocity field. Cotter et al.~\cite{cot1} illustrate their theoretical results by considering the velocity field for fluid flow generated by the two-dimensional Navier\textendash Stokes equation on a torus. They claim that {\em the case of Eulerian observations\textemdash direct observations of the velocity field itself\textemdash is then a model for weather forecasting} and that {\em the case of Lagrangian observations\textemdash observations of passive tracers advected by the flow\textemdash is then a model for data arising in oceanography.}

\section{The main results}
\label{mr}

We shall assume throughout the paper\textemdash unless stated otherwise\textemdash that $H, K$ are Hilbert spaces over the field ${\mathbb C}$ of complex numbers, that $(\Omega,\Sigma,\mu)$ is a probability space, and that $L^2(\Omega,H)$ and $L^2(\Omega,K)$ are the spaces of square-integrable random functions taking values in $H$ and $K$ respectively.

Let $\bff \in L^2(\Omega,H)$ and $\bfg \in L^2(\Omega,K)$ be zero-mean random functions.  We show that the auto-covariance $E_{\sbff \! \sbff} \in {\mathcal B}(H)$ is a nuclear operator.  If the range space $E_{\sbff\!\sbff}(H) \subseteq K$ is not closed we prove that the generalized inverse auto-covariance operator ${E_{\sbff \! \sbff}}^{\dag}:D({E_{\sbff \! \sbff}}^{\dag}) \subseteq H \rightarrow H$ is an unbounded, closed, densely defined, self-adjoint, linear operator.  We also show that the cross-covariance $E_{\sbfg \sbff} \in {\mathcal B}(H, K)$ is well defined and that the null space of $E_{\sbff \! \sbff}$ is a subspace of the null space of $E_{\sbfg \sbff}$.

Finally we show that there exists an optimal, closed, densely defined, linear operator $X = E_{\sbfg \sbff} {E_{\sbff \! \sbff}}^{\dag}:D({E_{\sbff\!\sbff}}^{\dag}) \subseteq H \rightarrow K$, a proximate observable function $\bfp = \bff_{\epsilon} \in L^2(\Omega,M)$ for each $\epsilon > 0$ with ${\mathbb E}[ \| \bfp - \bff \|^2] < \epsilon$ and $\bfp(\bomega) \in D({E_{\sbff\!\sbff}}^{\dag})$ for $\mu$-almost all $\bomega \in \Omega$, and a corresponding optimal estimate $\widehat{\bfg}_{\epsilon} = X\! \bfp \in L^2(\Omega,K)$ of the unobservable function $\bfg$ with mean square error
\begin{equation}
\label{minerror}
{\mathbb E}[ \| E_{\sbfg \sbff} {E_{\sbff \! \sbff}}^{\dag} \bfp - \bfg \|^2] = \mbox{tr}(E_{\sbfg \sbfg} - E_{\sbfg \sbfp}{E_{\sbfp \sbfp}}^{\dag} E_{\sbfp \sbfg}).
\end{equation}

The operator $X = E_{\sbfg \sbff}{E_{\sbff \! \sbff}}^{\dag}$ minimizes the mean square error $E[ \| X\! \bfp - \bfg \|^2]$  over all closed, densely defined, linear operators $X = T {E_{\sbff \! \sbff}}^{\dag}:D({E_{\sbff\!\sbff}}^{\dag}) \subseteq H \rightarrow K$ where $T \in {\mathcal B}(H,K)$.  The operator $X$ does not depend on the parameter $\epsilon$.  The notation $\bfp = \bff_{\epsilon}$ is simply a device to avoid the use of a double subscript in (\ref{minerror}). 

\section{Structure of the paper}
\label{sp}

In Section~\ref{pw} we review the previous work on this problem.  In Section~\ref{pre} we survey the necessary preliminary material.  We need to know that every Hilbert space has an orthonormal basis.  We state the relevant background theory \cite[pp 86\textendash 87]{yos1} and provide an example of a Hilbert space with an uncountable orthonormal basis.  In Section \ref{eno} we introduce an important elementary nuclear operator.   This material is taken from \cite{how2} but is central to later definitions and we need to repeat it here.  The necessary theory of the Bochner integral is summarized in Section \ref{birf}.  Once again we cite the text by Yosida~\cite[pp 130\textendash 134]{yos1}.   

The Hilbert space covariance operators are introduced and also justified in Section~\ref{co}.  We follow \cite{how2} but no longer assume that the Hilbert spaces are separable.  It is necessary to show that the auto-covariance is positive semi-definite and self-adjoint in order to extract a countable orthonormal basis for the orthogonal complement of the null space and thereby obtain an effective coordinate representation of the key operators.

The material in Section~\ref{giac} is new.  We show that the auto-covariance operator is nuclear and hence also compact.  We define the generalized inverse auto-covariance operator and show that in the general case it is an unbounded, closed, densely defined, self-adjoint, linear operator.  We also establish the standard properties of the generalized inverse auto-covariance operator and derive key formul{\ae} for the auto-covariance and cross-covariance of a specific linearly transformed random function that is used to establish the main result.  In Section~\ref{nsco} we show that the null space of the auto-covariance is a subspace of the null space of the cross-covariance.

In Section~\ref{sgep} we establish our main result\textemdash the solution to Problem~\ref{p1}.  The solution is presented in two parts.  Firstly we prove that a direct solution is possible if the observable function takes almost all values in the domain of the generalized inverse auto-covariance operator.  Secondly we argue that the direct solution is essentially preserved when the observable function is replaced by a suitable proximate observable function.  In Section~\ref{psp} we establish a key result, Lemma~\ref{gsrs}, that relates to practical aspects of the solution procedure.  To conclude, in Section~\ref{he}, we present a detailed study of a particular example.  The example highlights typical difficulties that arise when the results are applied.

\section{Previous work}
\label{pw}

Let $\bff \in L^2(\Omega, {\mathbb C}^m)$ and $\bfg \in L^2(\Omega, {\mathbb C}^n)$ be square-integrable, zero-mean, random vectors with realizations $\bff(\bomega) \in {\mathbb C}^m$ and $\bfg(\bomega) \in {\mathbb C}^n$ in finite-dimensional Euclidean space.  We assume that the covariance matrices
$$
E_{\sbff \! \sbff} = {\mathbb E}[\bff \bff^*] = \int_{\Omega} \bff(\bomega) {\bff(\bomega)}^* \mu(d\bomega) \in {\mathbb C}^{m \times m}
$$
and
$$
E_{\sbfg \sbff} = {\mathbb E}[\bfg \bff^*] = \int_{\Omega} \bfg(\bomega) {\bff(\bomega)}^* \mu(d\bomega) \in {\mathbb C}^{n \times m}
$$
are known, where ${\mathbb E}$ denotes the expectation operator.  If the matrix ${E_{\sbff \! \sbff}}^{-1}$ exists, then it has long been known \cite{sor1} that the best linear mean-square estimate $\widehat{\bfg} = X\!\bff$ of the random vector $\bfg$ from the observed data vector $\bff$ is
\begin{equation}
\label{oemat1}
\widehat{\bfg} = E_{\sbfg \sbff} {E_{\sbff \! \sbff}}^{-1} \bff
\end{equation}
with expected mean-square error
\begin{equation}
\label{omsemat1}
{\mathbb E}[ \| \widehat{\bfg} - \bfg \|^2] = \mbox{tr}( E_{\sbfg \sbfg} - E_{\sbfg \sbff}{E_{\sbff \! \sbff}}^{-1} E_{\sbff \sbfg})
\end{equation}
where $\mbox{tr}(\cdot)$ denotes the trace operator.  In this case the optimal solution $X = E_{\sbfg \sbff} {E_{\sbff \! \sbff}}^{-1} \in {\mathbb C}^{n \times m}$ is a finite-dimensional matrix and the linear mapping $\widehat{\bfg} = X\! \bff$ is defined by the relationship $\widehat{\bfg}(\bomega) = X\! \bff(\bomega) = X \cdot \bff(\bomega)$ for all $\bomega \in \Omega$.  Strictly speaking one should define an operator $L_X \in {\mathcal B}(L^2(\Omega,{\mathbb C}^m), L^2(\Omega,{\mathbb C}^n))$ by setting $[L_X \bff] (\bomega) = X \cdot \bff(\bomega)$ for each $\bomega \in \Omega$.  We prefer to write $X\! \bff$ rather than $L_X \bff$ so that $[X\! \bff](\bomega) = X \cdot \bff(\bomega)$ for all $\bomega \in \Omega$.  However we note that there are bounded linear transformations $F \in {\mathcal B}(L^2(\Omega,{\mathbb C}^m), L^2(\Omega,{\mathbb C}^n))$ that cannot be written in this way. 

Yamashita and Ogawa \cite{yam1} considered the special case $\bff = \bfg + \bfr$ where $\bff$ and $\bfr$ are independent random vectors with realizations in a finite-dimensional Euclidean space.  When the auto-covariance matrix $E_{\sbff \! \sbff}$ is singular they showed that an optimal estimate can be found in the form $\widehat{\bfg} = E_{\sbff \! \sbff} {E_{\sbff \! \sbff}}^{\dag}\bff$ where ${E_{\sbff \! \sbff}}^{\dag}$ is the Moore\textendash Penrose inverse \cite[Definition 2.2, p 10]{avr1}.  The expected mean-square error in this special case is ${\mathbb E}[ \| \widehat{\bfg} - \bfg \|^2] = {\mathbb E}[ \| \bfr \|^2] = \mbox{tr}(E_{\sbfr \sbfr})$.   Hua and Liu \cite{hua1} improved this result by showing that the random vectors $\bff$ and $\bfg$ can lie in different spaces and that no special relationship between the two vectors is necessary.  The optimal estimate is now given by
\begin{equation}
\label{oemat2}
\widehat{\bfg} = E_{\sbfg \sbff} {E_{\sbff \! \sbff}}^{\dag} \bff
\end{equation}
with expected mean-square error
\begin{equation}
\label{omsemat2}
{\mathbb E}[ \| \widehat{\bfg} - \bfg \|^2] =  \mbox{tr}(E_{\sbfg \sbfg} - E_{\sbfg \sbff} {E_{\sbff \! \sbff}}^{\dag} E_{\sbff \sbfg}).
\end{equation}
This solution was extended to random vectors taking values in different Hilbert spaces by Fomin and Ruzhansky \cite[Theorem 4.1]{fom1} and by Howlett, Pearce and Torokhti \cite[Theorem 3]{how2}, independently, and at about the same time.  In each case the authors assumed that the generalized inverse auto-covariance operator ${E_{\sbff \! \sbff}}^{\dag}$ was a bounded linear operator.  We make no such assumption here and propose a more general solution procedure that allows the generalized inverse operator ${E_{\sbff \!\sbff}}^{\dag}$ to be unbounded.  This relaxation has profound implications.  See our earlier remarks in Sections \ref{intro} and \ref{mr}.

\section{Preliminaries}
\label{pre}

A substantial portion of the preliminary material in Sections \ref{eno} and \ref{birf} is reprised from \cite{how2}.  We begin with some basic facts about Hilbert space.  In particular we need to know that every Hilbert space has an orthonormal basis which may or may not be countable.  We follow the presentation in Yosida \cite[pp 86\textendash 87]{yos1}.

\begin{defn}
\label{def1yos}
A set $S$ of vectors in a Hilbert space $H$ is called an orthogonal set if $\langle \bfx, \bfu \rangle = 0$ for all $\bfx, \bfu \in S$ with $\bfx \neq \bfu$.  If, in addition, $\| \bfx \| = 1$ for all $\bfx \in S$ then we say the $S$ is an orthonormal set.  An orthonormal set $S$ of a Hilbert space $H$ is called a complete orthonormal system or an orthonormal basis of $H$, if no orthonormal set of $H$ contains $S$ as a proper subset.  $\hfill \Box$
\end{defn}

Some authors say that a {\em complete} orthonormal set is a {\em maximal} orthonormal set.  See Naylor and Sell \cite[Definition 5.17.4, p 306]{nay1}. 

\begin{thm}
\label{t1yos}
A Hilbert space $H$ containing a non-zero vector has at least one complete orthonormal system.  Moreover, if $S$ is any orthonormal set in $H$, there is a complete orthonormal set containing $S$. $\hfill \Box$
\end{thm}

\begin{thm}
\label{t2yos}
Let $S = \{ \bfx_{\alpha}\}_{\alpha \in A}$ be a complete orthonormal system of a Hilbert space $H$.  For any $\bfh \in H$ we define the Fourier coefficients of $\bfh$ with respect to $S$ by $h_{\alpha} = \langle \bfh, \bfx_{\alpha} \rangle$ for each $\alpha \in A$.  Then we have Parseval's relation $\| \bfh \|^2 = \sum_{\alpha \in A} |h_{\alpha}|^2$. $\hfill \Box$.
\end{thm}

\begin{cor}
\label{cor1yos}
Let $S = \{ \bfx_{\alpha} \}_{\alpha \in A}$ be a complete orthonormal system in $H$.  For each $\bfh \in H$ there is a countable subset $S_{\sbfh,+} \subseteq S$ such that $h_{\alpha} = \langle \bfh, \bfx_{\alpha} \rangle \neq 0$ for $\alpha \in S_{\sbfh,+}$ and $h_{\alpha} = \langle \bfh, \bfx_{\alpha} \rangle = 0$ for $\alpha \in S_{\sbfh,0} = S \setminus S_{\sbfh,+}$.  If we write $S_{\sbfh,+}$ in the form $S_{\sbfh,+} = \{ \bfx_{\sbfh,j} \}_{j \in {\mathbb N}}$ for convenience then we have $\| \sum_{j = n+1}^{\infty}  \langle \bfh, \bfx_{\sbfh, j} \rangle \bfx_{\sbfh,j} \|^2 \rightarrow 0$ as $n \rightarrow \infty$ and we can represent $\bfh$ by the Fourier series $\bfh = \sum_{j \in {\mathbb N}}  \langle \bfh, \bfx_{\sbfh, j} \rangle \bfx_{\sbfh,j}$.  $\hfill \Box$
\end{cor}

The following example is taken from Naylor and Sell \cite[Example 10, p 320]{nay1}.
   
\begin{exa}
\label{ex1}
The set $AP$ of all complex-valued almost periodic functions $\bff:{\mathbb R} \rightarrow {\mathbb C}$ with the property
$$
\lim_{T \rightarrow \infty} (1/T) \mbox{$\int_{[-T,T]}$} |\bff(t)|^2 dt < \infty
$$
becomes a Hilbert space if we define an inner product
$$
\langle \bff, \bfg \rangle = \lim_{T \rightarrow \infty} (1/T) \mbox{$\int_{[-T,T]}$} \bff(t) \overline{\bfg(t)} dt
$$
for each $\bff, \bfg \in AP$ and an associated norm $\| \bff \| = \langle \bff, \bff \rangle^{1/2}$ for each $\bff \in AP$.  The set $\{ \bfe_{\alpha} \}_{\alpha \in {\mathbb R}}$ defined by $\bfe_{\alpha}(t) = e^{i \alpha t}$ for each $t \in {\mathbb R}$ forms an uncountable orthonormal basis for $AP$. $\hfill \Box$
\end{exa}  

\subsection{An elementary nuclear operator.}
\label{eno}

For each $\bfh \in H$ define a corresponding linear operator $J_{\sbfh} \in {\mathcal B}({\mathbb C}, H)$ by the formula $J_{\sbfh} z = z \bfh$.  The range space $J_h({\mathbb C}) \subseteq H$ is a one-dimensional subspace spanned by $\bfh$.  The adjoint operator ${J_{\sbfh}}^* \in {\mathcal B}(H, {\mathbb C})$ is defined by the relationship
$$
 \overline{z} {J_{\sbfh}}^*\bfx = \langle {J_{\sbfh}}^*\bfx, z \rangle = \langle \bfx, J_{\sbfh}z \rangle = \langle \bfx, z \bfh \rangle = \overline{z} \langle \bfx, \bfh \rangle 
$$
for all $\bfx \in H$ and $z \in {\mathbb C}$ and hence ${J_{\sbfh}}^*\bfx = \langle \bfx, \bfh \rangle$ for each $\bfx \in H$.  If $\bfx \perp \bfh$ then ${J_{\sbfh}}^*\bfx = 0$.   If $T \in {\mathcal B}(H,K)$ and we define $\bfk = T \bfh$ then $J_{\sbfk} \in {\mathcal B}({\mathbb C}, K)$ and we have $J_{\sbfk}z = z \bfk = z T \bfh = T (z \bfh) = TJ_{\sbfh}z$ for all $z \in {\mathbb C}$.  Thus $J_{\sbfk} = TJ_{\sbfh}$.  We also have ${J_{\sbfk}}^* = {J_{\sbfh}}^* T^* \in {\mathcal B}(K, {\mathbb C})$ and $J_{\sbfk}{J_{\sbfk}}^* = T J_{\sbfh}{J_{\sbfh}}^* T^* \in {\mathcal B}(K)$.  If $\bfh \in H$ and $\bfk \in K$ the operator $J_{\sbfk}{J_{\sbfh}}^* \in {\mathcal B}(H, J_{\sbfk}({\mathbb C}))$ is given by $J_{\sbfk}{J_{\sbfh}}^*\bfx = \langle \bfx, \bfh \rangle \bfk$ for each $\bfx \in H$ and so $\langle J_{\sbfk}{J_{\sbfh}}^*\bfx, \bfy \rangle = \langle \bfx, \bfh \rangle  \langle \bfk, \bfy, \rangle$ for each $\bfx \in H$ and $\bfy \in K$.

We are particularly interested in the operator $J_{\sbfh}{J_{\sbfh}}^* \in {\mathcal B}(H, J_{\sbfh}({\mathbb C}))$.  Since $J_{\sbfh}({\mathbb C}) \subseteq H$ is a one-dimensional subspace it follows that $J_{\sbfh}{J_{\sbfh}}^*$ is a compact operator \cite[pp 379\textendash 381]{nay1}.  If $\bfx \in J_{\sbfh}({\mathbb C})$ then $\bfx = w \bfh$ for some $w \in {\mathbb C}$ and so $J_{\sbfh} {J_{\sbfh}}^* \bfx = \langle w \bfh, \bfh \rangle \bfh = \| \bfh \|^2 w \bfh = \| \bfh \|^2 \bfx$.  Thus $\bfx$ is an eigenvector with corresponding eigenvalue $\| \bfh\|^2$.  If $\bfu \in J_{\sbfh}({\mathbb C})^{\perp}$ then $J_{\sbfh} {J_{\sbfh}}^*\bfu = \langle \bfu, \bfh \rangle \bfh = \bfzero$ and so $\bfu$ is an eigenvector with corresponding eigenvalue $0$.  Write $H = J_{\sbfh}({\mathbb C}) \oplus J_{\sbfh}({\mathbb C})^{\perp}$.  Define $\bfx = \bfh/\|\bfh\|$ and let $\{\bfu_{\alpha} \}_{\alpha \in A_0}$ be a complete orthonormal set in $J_{\sbfh}({\mathbb C})^{\perp}$.  The trace of the positive semi-definite, self-adjoint operator $J_{\sbfh}{J_{\sbfh}}^* \in {\mathcal B}(H)$ is given by
\begin{eqnarray*}
\lefteqn{\mbox{tr}( J_{\sbfh}{J_{\sbfh}}^*) = \langle J_{\sbfh}{J_{\sbfh}}^* \bfx, \bfx \rangle + \mbox{$\sum_{\alpha \in A_0}$} \langle J_{\sbfh}{J_{\sbfh}}^* \bfu_{\alpha}, \bfu_{\alpha} \rangle } \hspace{4cm} \\
& & = \langle \bfh, \bfx \rangle \langle \bfx, \bfh \rangle = \langle \bfh, \bfh \rangle^2/\| \bfh \|^2 = \| \bfh \|^2 < \infty.
\end{eqnarray*}
Thus $J_{\sbfh}{J_{\sbfh}}^*$ is a {\em nuclear} or equivalently {\em trace-class} operator \cite{bal1, dun1, yos1}.

\subsection{The Bochner integral of a random function.}
\label{birf}

Let $X$ be a Banach space over the field ${\mathbb C}$ of complex numbers with norm $\| \cdot \|: X \rightarrow [0, \infty)$.  We say that a function $\bff : \Omega \rightarrow X$ is a vector-valued random function or simply a  random function.  The following definitions and results have been extracted from the text by Yosida \cite[pp 130\textendash134]{yos1}.
 
\begin{defn}
\label{def2yos}
The random function $\bff:\Omega \rightarrow X$ is said to be finitely valued if there exists a finite collection of disjoint sets $\{E_j\}_{j=1}^m \in \Sigma$ such that $\bff(\bomega) = \bfc_j$ for each $\bomega \in E_j$ and each $j=1,2,\ldots,m$ and $\bff(\bomega) = \bfzero$ elsewhere.  In such cases we define the $\mu$-integral of $\bff$ by the formula $\int_{\Omega} \bff(\bomega) \mu (d \bomega) = \sum_{j=1}^m \bfc_j \mu(E_j)$.   $\hfill \Box$
\end{defn}

\begin{defn}
\label{def3yos}
The function $\bff: \Omega \rightarrow X$ is strongly $\Sigma$-measurable if there exists a sequence $\{ \bff_n \}_{n \in {\mathbb N}}$ of finitely-valued functions $\bff_n:\Omega \rightarrow H$ with $\| \bff(\bomega) - \bff_n(\bomega) \| \rightarrow 0$ for $\mu$-almost all $\bomega \in \Omega$. $\hfill \Box$
\end{defn}

\begin{defn}
\label{def4yos}
The function $\bff: \Omega \rightarrow X$ is Bochner $\mu$-integrable if there exists a sequence $\{ \bff_n\}_{n \in {\mathbb N}}$ of finitely-valued functions $\bff_n: \Omega \rightarrow X$ with $\| \bff_n(\bomega) - \bff(\bomega) \| \rightarrow 0$ for $\mu$-almost all $\bomega \in \Omega$ in such a way that
$$
\lim_{n \rightarrow \infty} \mbox{$\int_{\Omega}$}\, \| \bff_n(\bomega) - \bff(\bomega) \| \mu(d \bomega) = 0.
$$
For each set $E \in \Sigma$ the Bochner $\mu$-integral of $\bff(\bomega)$ over $S$ is defined by
$$
\mbox{$\int_E \bff(\bomega)$} \mu(d \bomega) = \lim_{n \rightarrow \infty} \mbox{$\int_{\Omega}$}\, \chi_E(\bomega) \bff_n(\bomega) \mu(d \bomega)
$$
where $\chi_E:\Omega \rightarrow \{0,1\}$ is the characteristic function for $E$ given by $\chi_E(\bomega) = 1$ for $\bomega \in E$ and $\chi_E(\bomega) = 0$ otherwise. $\hfill \Box$
\end{defn} 

\begin{thm}
\label{t3yos}
A strongly $\Sigma$-measurable function $\bff:\Omega \rightarrow X$ is Bochner $\mu$-integrable if and only if the function $\| \bff \|: \Omega \rightarrow [0, \infty)$ defined by $\| \bff \|(\bomega) = \| \bff(\bomega) \|$ for all $\bomega \in \Omega$ is $\mu$-integrable in which case
$$
\| \mbox{$\int_E$}\, \bff(\bomega) \mu(d \bomega)\ \| \leq \mbox{$\int_E$}\, \|\bff(\bomega)\| \mu(d \bomega)
$$
for each $E \in \Sigma$. $\hfill \Box$
\end{thm}

\begin{cor}
\label{cor2yos}
Let $X$ and $Y$ be Banach spaces and suppose that $T \in {\mathcal B}(X,Y)$.  If the function $\bff:\Omega \rightarrow X$ is Bochner $\mu$-integrable then the function $\bfg = T\bff:\Omega \rightarrow Y$ defined by $\bfg(\bomega) = T \bff(\bomega)$ for $\mu$-almost all $\bomega \in \Omega$ is Bochner $\mu$-integrable with
$$
\mbox{$\int_E$}\, \bfg(\bomega) \mu(d \bomega) = T \mbox{$\int_E$}\, \bff(\bomega) \mu(d \bomega)
$$
for each $E \in \Sigma$.  $\hfill \Box$
\end{cor}

Let $\bff:\Omega \rightarrow X$ be a Bochner $\mu$-integrable random function taking values in the Banach space $X$. The expected value of $\bff$ is defined by
$$
{\mathbb E}[\bff] = \mbox{$\int_{\Omega}$}\, \bff(\bomega) \mu(d\bomega)
$$
and we note from Theorem \ref{t3yos} that $\|{\mathbb E}[\bff]\| \leq {\mathbb E}[\| \bff \|]$. When $T \in {\mathcal B}(X,Y)$ is a bounded linear map from the Banach space $X$ to the Banach space $Y$, it follows from Corollary \ref{cor2yos} that ${\mathbb E}[T \bff] = T{\mathbb E}[\bff]$.

The theory of random functions in Hilbert space is an extension of the corresponding theory in Banach space. Of particular interest are those properties relating to the scalar product which are used directly in defining the special operators for the optimal filter.  Let $H$ be a Hilbert space with scalar product $\langle \cdot,\cdot \rangle$ and let $\bff:\Omega \rightarrow H$ be a finitely-valued random function defined by $\bff(\bomega) = \sum_{j=1}^m \chi_j(\bomega) \bfc_j$ where $\{E_j\}_{j=1}^m$ are disjoint $\mu$-measurable sets and $\chi_j:\Omega \rightarrow \{0, 1\}$ is the characteristic function for $E_j$ for each $j=1,\ldots,m$.  Since $\| u(\bomega )\| ^2 = \mbox{$\sum _{j=1}^m$}\, \chi _j(\bomega )\| \bfc_j \| ^2$, it follows that if $T \in {\mathcal B}(H)$ is a bounded linear map, then we can use the elementary inequalities
$$
|\langle \bfc_{j},T[\bfc_{k}] \rangle| \leq \|T\| \cdot \|\bfc_{j}\| \cdot \|\bfc_{k}\| \quad \mbox{and} \quad \|\bfc_{j}\| \cdot \|\bfc_{k}\| \leq \left[ \|\bfc_{j}\|^{2} + \|\bfc_{k}\|^{2} \right]/2
$$
to deduce that
\begin{eqnarray*}
\lefteqn{\langle \mbox{$\int _{\Omega}$}\, \bff(\bomega)\mu (d\bomega ), \mbox{$\int _{\Omega}$} T[\bff(\bomega)]\mu (d\bomega ) \rangle } \hspace{3cm} \\
& = & \mbox{$\sum _{j=1}^m\sum _{k=1}^m$}\, \mu (E_j)\mu (E_k) \langle \bfc _j,T[\bfc _k]\rangle \\
& = & \| T \| \mbox{$\sum _{j=1}^m\sum _{k=1}^m$}\, \mu (E_j)\mu (E_k) \cdot \| \bfc _j\| \cdot \| \bfc _k\| \\
& \leq & \| T \| \mbox{$\sum _{j=1}^m\sum _{k=1}^m$}\, \mu (E_j)\mu (E_k)\| \cdot (\| \bfc_j\|^2 + \| \bfc _k\|^2)/2 \\
& = & \| T \| \mbox{$\sum _{j=1}^m$}\, \mu (E_j) \| \bfc _j\| ^2 \\
& = & \| T \| \mbox{$\int _{\Omega}$}\, \| \bff(\bomega )\| ^2\mu (d\bomega ).
\end{eqnarray*}
By taking appropriate limits, we can extend the above argument to establish the following general results, which are used to justify construction of the optimal filter.

\begin{thm}
\label{squareint1}
Let $H$ be a Hilbert space.  If the random function $\bff:\Omega \rightarrow H$ is strongly $\Sigma$\textendash measurable and $\|\bff\|^{2}:\Omega \rightarrow [0,\infty)$ is $\mu$-integrable, then $\bff$ is Bochner $\mu$-integrable and for each bounded linear map $T \in {\mathcal B}(H)$ we have
$$
\langle \mbox{$\int _{\Omega}$}\, \bff(\bomega )\mu (d\bomega ), \mbox{$\int _{\Omega}$} T\bff(\bomega) \mu (d\bomega ) \rangle \leq \| T \| \mbox{$\int _{\Omega}$}\, \| \bff(\bomega)\| ^2\mu (d\bomega ).
$$
$\hfill \Box$
\end{thm}

\begin{cor}
\label{squareint2}
If $\bff:\Omega \rightarrow H$ is strongly $\Sigma$-measurable and $\|\bff\|^{2}:\Omega \rightarrow [0,\infty)$ is $\mu$-integrable, then
$$
\| \mbox{$\int _{\Omega}$}\, \bff(\bomega) \mu (d\bomega) \|^2 \leq \mbox{$\int_{\Omega}$}\, \| \bff(\bomega ) \|^2 \mu (d\bomega).
$$
$\hfill \Box$
\end{cor}

Theorem \ref{squareint1} and Corollary \ref{squareint2} can be expressed in terms of expected values. Let $T \in {\mathcal B}(H)$ and let $\bff:\Omega \rightarrow X$ be a random function. If $\|\bff\|^{2}: \Omega \rightarrow [0, \infty)$ is $\mu$-integrable, then $\langle {\mathbb E}[\bff], {\mathbb E}[T \bff] \rangle \leq \|T\| \cdot {\mathbb E}[\|\bff\|^{2}]$ and $\|{\mathbb E}[\bff]\|^{2} \leq {\mathbb E}[\|\bff\|^{2}]$.

\section{The covariance operators}
\label{co}

If $\bff:\Omega \rightarrow H$ is strongly $\Sigma$-measurable and ${\mathbb E}[ \| \bff\|^2] < \infty$ then we say that $\bff(\bomega)$ is $\mu$-square-integrable on $\Omega$ and we write $\bff \in L^2(\Omega,H)$.   If $\bff_1, \bff_2 \in L^2(\Omega,H)$ and we define the inner product $\dlangle \bff_1, \bff_2 \drangle = {\mathbb E}[ \langle \bff_1, \bff_2 \rangle]$ then $L^2(\Omega,H)$ becomes a Hilbert space.  For each $\bff \in L^2(\Omega, H)$ we write $\dnorm \bff \dnorm = \dlangle \bff, \bff \drangle^{1/2} = {\mathbb E}[ \| \bff\|^2]^{1/2}$ for the corresponding norm.  If $ \bfx \in H$ and we define an associated constant function $\bfx:\Omega \rightarrow H$ by setting $\bfx(\bomega) = \bfx$ for all $\bomega \in \Omega$ then $\dnorm \bfx \dnorm^2 = {\mathbb E}[\| \bfx \|^2] = \int_{\Omega} \| \bfx \|^2 \mu(d \bomega) = \| \bfx \|^2$.  Thus $\bfx \in L^2(\Omega,H)$.  Similarly if $\bfx, \bfu \in H$ then $\dlangle \bfx, \bfu \drangle = \int_{\Omega} \langle \bfx, \bfu \rangle \mu(d \bomega) = \langle \bfx, \bfu \rangle$.  Thus we could regard $H$ as a subspace of $L^2(\Omega, H)$. 

\subsection{The basic pointwise functions.}
\label{bpf}

Suppose that $\bff \in L^2(\Omega,H)$ is a random function with zero mean.  For each $\bomega \in \Omega$ we have $J_{\sbff(\sbfomega)} \in {\mathcal B}({\mathbb C},H)$ defined by $J_{\sbff(\sbfomega)}z = z \bff(\bomega)$ for all $z \in {\mathbb C}$ and ${J_{\sbff(\sbfomega)}}^* \in {\mathcal B}(H, {\mathbb C})$ defined by ${J_{\sbff(\sbfomega)}}^*\bfx = \langle \bfx, \bff(\bomega) \rangle$ for each $\bfx \in H$.  Therefore $J_{\sbff(\sbfomega)}{J_{\sbff(\sbfomega)}}^* \in {\mathcal B}(H)$ for all $\bomega \in \Omega$ with $J_{\sbff(\sbfomega)}{J_{\sbff(\sbfomega)}}^*\bfx = \langle \bfx, \bff(\bomega) \rangle \bff(\bomega)$ and 
$$
\langle J_{\sbff(\sbfomega)} {J_{\sbff(\sbfomega)}}^* \bfx, \bfu \rangle = \langle \bfx, \bff(\bomega) \rangle \langle \bff(\bomega), \bfu \rangle
$$
for all $\bomega \in \Omega$ and each $\bfx, \bfu \in H$.  We also have ${J_{\sbff(\sbfomega)}}^* J_{\sbff(\sbfomega)} \in {\mathcal B}({\mathbb C})$ for all $\bomega \in \Omega$ with ${J_{\sbff(\sbfomega)}}^* J_{\sbff(\sbfomega)} z = z \| \bff(\bomega)\|^2$ for all $\bomega \in \Omega$ and each $z \in {\mathbb C}$.  If $T \in {\mathcal B}(H, K)$ then $J_{T \sbff(\sbfomega)} = T J_{\sbff(\sbfomega)}$ and ${J_{T \sbff(\sbfomega)}}^* = {J_{\sbff(\sbfomega)}}^* T^*$ for all $\bomega \in \Omega$.  To continue we must show that certain key functions are measurable. 

\begin{lem}
\label{pmeas}
Let $\bfx \in H$ and $\bff \in L^2(\Omega,H)$.  If we define an associated random function $\bfp:\Omega \rightarrow H$ by setting $\bfp(\bomega) = J_{\sbff(\sbfomega)}{J_{\sbff(\sbfomega)}}^*\bfx = \langle \bfx, \bff(\bomega) \rangle \bff(\bomega)$ for all $\bomega \in \Omega$ then $\bfp$ is strongly $\Sigma$-measurable. $\hfill \Box$
\end{lem}

{\bf Proof} \ Let $\{ \bff_n \}_{n \in {\mathbb N}}$ be a sequence of finitely valued functions such that $\| \bff_n(\bomega) - \bff(\bomega) \| \rightarrow 0$ as $n \rightarrow \infty$ for $\mu$-almost all $\bomega \in \Omega$.  Define $\bfp_n:\Omega \rightarrow H$ by setting $\bfp_n(\bomega) = \langle \bfx, \bff_n(\bomega) \rangle \bff_n(\bomega)$ for all $\bomega \in \Omega$.  Then $\{ \bfp_n \}_{n \in {\mathbb N}}$ is a sequence of finitely valued functions with
\begin{eqnarray*}
\| \bfp_n(\bomega) - \bfp(\bomega) \| & = & \| \langle \bfx, \bff_n(\bomega) \rangle \bff_n(\bomega) - \langle \bfx, \bff(\bomega) \rangle \bff(\bomega) \| \\
& & \hspace{-2.5cm} = \| \langle \bfx, \bff_n(\bomega) - \bff(\bomega) \rangle \bff_n(\bomega) + \langle \bfx, \bff(\bomega) \rangle [\bff_n(\bomega) - \bff(\bomega)] \| \\
& & \hspace{-2.5cm} \leq \| \bfx \| \cdot \| \bff_n(\bomega) - \bff(\bomega)\| \cdot \| \bff_n \| + \| \bfx \| \cdot \| \bff(\bomega) \| \cdot \| \bff_n(\bomega) - \bff(\bomega) \| \rightarrow 0
\end{eqnarray*}
as $n \rightarrow \infty$ for $\mu$-almost all $\bomega \in \Omega$.  Therefore $\bfp$ is strongly $\Sigma$-measurable.  $\hfill \Box$

\subsection{The auto-covariance operator.}
\label{aco}

Suppose that $\bff \in L^2(\Omega,H)$ is a $\mu$-square-integrable random function with zero mean.  The inequality
$$
\| \mbox{$\int_{\Omega}$} \langle \bfx, \bff(\bomega) \rangle \bff(\bomega) \mu(d \bomega)  \| \leq \| \bfx \| \mbox{$\int_{\Omega}$} \| \bff(\bomega) \|^2 \mu(d \bomega) = \| \bfx \| \cdot \dnorm \bff \dnorm^2 < \infty
$$
justifies the definition of an operator $E_{\sbff \! \sbff} \in {\mathcal B}(H)$ by setting
$$
E_{\sbff \! \sbff} \bfx = \mbox{$\int_{\Omega}$} J_{\sbff(\sbfomega)}{J_{\sbff(\sbfomega)}}^* \bfx\, \mu(d \bomega) = \mbox{$\int_{\Omega}$} \langle \bfx, \bff(\bomega) \rangle \bff(\bomega)\, \mu(d \bomega) = {\mathbb E}[ \langle \bfx, \bff \rangle \bff]
$$
for all $\bfx \in H$.  Let $T \in {\mathcal B}(H,K)$.  We have
\begin{eqnarray*}
T E_{\sbff \! \sbff}\, T^* \bfy & = & T \mbox{$\int_{\Omega}$} \langle T^*\bfy, \bff(\bomega) \rangle \bff(\bomega)\, \mu(d \bomega) \\
& = & \mbox{$\int_{\Omega}$} \langle T^*\bfy, \bff(\bomega) \rangle T \cdot \bff(\bomega)\, \mu(d \bomega) \\
& = & \mbox{$\int_{\Omega}$} \langle \bfy, T\! \bff(\bomega) \rangle T\! \bff (\bomega)\, \mu(d \bomega) \\
& = & E_{T \! \sbff, T \! \sbff}\, \bfy
\end{eqnarray*}
for all $\bfy \in K$.  Thus we have $E_{T \! \sbff, T \! \sbff} = T E_{\sbff \! \sbff}\, T^* \in {\mathcal B}(K)$.  We also have
\begin{eqnarray*}
\langle E_{\sbff \! \sbff} \bfx, \bfu \rangle & = & \mbox{$\int_{\Omega}$} \langle J_{\sbff(\sbfomega)} {J_{\sbff(\sbfomega)}}^* \bfx, \bfu \rangle \mu(d \bomega) \\
& = & \mbox{$\int_{\Omega}$} \langle \bfx, \bff(\bomega) \rangle \langle \bff(\bomega), \bfu \rangle \mu(d \bomega) = {\mathbb E}[\langle \bfx, \bff \rangle \langle \bff, \bfu \rangle]
\end{eqnarray*}
for all $\bfx, \bfu \in H$.  Therefore
\begin{eqnarray*}
\langle E_{\sbff \! \sbff} \bfx, \bfx \rangle & = & \mbox{$\int_{\Omega}$} \langle J_{\sbff(\sbfomega)} {J_{\sbff(\sbfomega)}}^* \bfx, \bfx \rangle \mu(d \bomega) \\
& = & \mbox{$\int_{\Omega}$} | \langle \bfx, \bff(\bomega) \rangle |^2 \mu(d \bomega) = {\mathbb E}[| \langle \bfx, \bff \rangle |^2] \geq 0
\end{eqnarray*}
and hence $E_{\sbff \! \sbff}$ is positive semi-definite and self-adjoint.  We have the following elementary, but important, results.

\begin{lem}
\label{aclem1}
Let $\bff \in L^2(\Omega,H)$ and let $\bfx \in H$.  Then $\bfx \in {E_{\sbff \! \sbff}}^{-1}( \{ \bfzero \})$ if and only if $ \langle \bfx, \bff(\bomega) \rangle = 0$ for $\mu$-almost all $\bomega \in \Omega$. $\hfill \Box$
\end{lem}

{\bf Proof}\quad If $\bfx \in {E_{\sbff \! \sbff}}^{-1}( \{ \bfzero \})$ then $E_{\sbff \! \sbff} \bfx = \bfzero$ and so
$$
\langle E_{\sbff \! \sbff} \bfx, \bfx \rangle = \mbox{$\int_{\Omega}$} | \langle \bfx, \bff(\bomega) \rangle |^2 \mu (d \bomega) = 0.
$$
Therefore $\langle \bfx, \bff(\bomega) \rangle = 0$ for $\mu$-almost all $\bomega \in \Omega$.  Conversely if $\langle \bfx, \bff(\bomega) \rangle = 0$ for $\mu$-almost all $\bomega \in \Omega$ then
$$
\langle E_{\sbff \! \sbff} \bfx, \bfu \rangle = \mbox{$\int_{\Omega}$} \langle \bfx, \bff(\bomega) \rangle \langle \bff(\bomega), \bfu \rangle \mu (d \bomega) = 0
$$
for all $\bfu \in H$.  Therefore $E_{\sbff \! \sbff} \bfx = \bfzero$ and hence $\bfx \in {E_{\sbff \! \sbff}}^{-1}( \{ \bfzero \})$.  $\hfill \Box$

\begin{lem}
\label{aclem2}
Let $\bff \in L^2(\Omega,H)$.  Then
$$
\mbox{\rm tr}(E_{\sbff \! \sbff}) = \mbox{$\int_{\Omega}$} \mbox{\rm tr}( J_{\sbff(\sbfomega)} {J_{\sbff(\sbfomega)}}^*)\, \mu(d \bomega).
$$
$\hfill \Box$
\end{lem}

{\bf Proof} \quad Let $\{\bfe_{\beta} \}_{\beta \in B}$ be a complete orthonormal set in $H$.  Since 
$$
\mbox{tr}(E_{\sbff \! \sbff}) = \mbox{$\sum_{\beta \in B}$} \langle E_{\sbff \! \sbff}\, \bfe_{\beta}, \bfe_{\beta} \rangle = {\mathbb E}[ \| \bff \|^2] = \dnorm \bff \dnorm^2 < \infty
$$
there is at most a countable subset $B_+ \subseteq B$ with $\langle E_{\sbff \! \sbff} \bfe_{\beta}, \bfe_{\beta} \rangle > 0$ for each $\beta \in B_+$.  Lemma \ref{aclem1} shows that $\langle E_{\sbff \! \sbff} \bfe_{\beta}, \bfe_{\beta} \rangle = 0$ if and only if $\langle J_{\sbff(\sbfomega)}{J_{\sbff(\sbfomega)}}^* \bfe_{\beta}, \bfe_{\beta} \rangle = | \langle \bfe_{\beta}, \bff(\bomega) \rangle |^2 = 0$ for $\mu$-almost all $\bomega \in \Omega$ in which case $\beta \in B_0 = B \setminus B_+$. Now we have
\begin{eqnarray*}
\mbox{tr}(E_{\sbff \! \sbff}) & = & \mbox{$\sum_{\beta \in B_+}$} \langle E_{\sbff \! \sbff} \bfe_{\beta}, \bfe_{\beta} \rangle \\
& = & \mbox{$\sum_{\beta \in B_+} \int_{\Omega}$} \langle J_{\sbff(\sbfomega)}{J_{\sbff(\sbfomega)}}^* \bfe_{\beta}, \bfe_{\beta} \rangle \mu(d \bomega) \\
& = & \mbox{$\int_{\Omega} \sum_{\beta \in B_+}$} \langle J_{\sbff(\sbfomega)}{J_{\sbff(\sbfomega)}}^* \bfe_{\beta}, \bfe_{\beta} \rangle \mu(d \bomega) \\
& = & \mbox{$\int_{\Omega}$} \mbox{tr}( J_{\sbff(\sbfomega)}{J_{\sbff(\sbfomega)}}^*) \mu(d \bomega)
\end{eqnarray*}
as required. $\hfill \Box$

\subsection{The cross-covariance operator.}
\label{cco}

Suppose that $\bff \in L^2(\Omega,H)$ and $\bfg \in L^2(\Omega,K)$ are $\mu$-square-integrable random functions with zero mean.  By essentially repeating previous arguments we deduce that $J_{\sbfg(\sbfomega)}{J_{\sbff(\sbfomega)}}^* \in {\mathcal B}(H, K)$ with $J_{\sbfg(\sbfomega)}{J_{\sbff(\sbfomega)}}^*\bfx = \langle \bfx, \bff(\bomega) \rangle \bfg(\bomega) \in K$ for all $\bomega \in \Omega$ and each $\bfx \in H$.  It follows that for fixed $\bfx \in H$ the function $\bfq:\Omega \rightarrow K$ defined by $\bfq(\bomega) = J_{\sbfg(\sbfomega)}{J_{\sbff(\sbfomega)}}^*\bfx = \langle \bfx, \bff(\bomega) \rangle \bfg(\bomega)$ for all $\bomega \in \Omega$ is strongly $\Sigma$-measurable.  Now the inequality
$$
\| \mbox{$\int_{\Omega}$} \langle \bfx, \bff(\bomega) \rangle \bfg(\bomega) \mu(d \bomega) \| \leq \|\bfx\| \cdot \dnorm \bff \dnorm \cdot \dnorm \bfg \dnorm  < \infty
$$
justifies the definition of an operator $E_{\sbfg \sbff} \in {\mathcal B}(H,K)$ by the formula
$$
E_{\sbfg \sbff}\bfx = \mbox{$\int_{\Omega}$} J_{\sbfg(\sbfomega)}{J_{\sbff(\sbfomega)}}^* \bfx\, \mu(d \bomega) = \mbox{$\int_{\Omega}$} \langle \bfx, \bff(\bomega) \rangle \bfg(\bomega)\, \mu(d \bomega) = {\mathbb E}[ \langle \bfx, \bff \rangle \bfg]
$$
for each $\bfx \in H$. We also have $\langle J_{\sbfg(\sbfomega)}{J_{\sbff(\sbfomega)}}^*\bfx, \bfy \rangle = \langle \bfx, \bff(\bomega) \rangle \langle \bfg(\bomega), \bfy \rangle$ for all $\bomega \in \Omega$ and each $\bfx \in H$ and $\bfy \in K$ and so
\begin{eqnarray*}
\langle E_{\sbfg \sbff}\bfx, \bfy \rangle & = & \mbox{$\int_{\Omega}$} \langle J_{\sbfg(\sbfomega)}{J_{\sbff(\sbfomega)}}^*\bfx, \bfy \rangle \mu(d \bomega) \\
& = & \mbox{$\int_{\Omega}$} \langle \bfx, \bff(\bomega) \rangle \langle \bfg(\bomega), \bfy \rangle \mu(d \bomega) = {\mathbb E}[ \langle \bfx, \bff \rangle \langle \bfg, \bfy \rangle]
\end{eqnarray*}
for each $\bfx \in H$ and $\bfy \in K$.  If $\bfg, \bfk \in L^2(\Omega,K)$ we can use the definitions and basic algebra to show that $E_{\sbfg + \sbfk,\, \sbfg + \sbfk} = E_{\sbfg \sbfg} + E_{\sbfk \sbfg} + E_{\sbfg \sbfk} + E_{\sbfk \sbfk}$. 

\subsection{The definitive properties of the covariance operators.}
\label{dpco}

The operator $E_{\sbff \! \sbff} \in {\mathcal B}(H)$ is self-adjoint and positive semi-definite.  Thus we can find a countable  orthonormal basis of eigenvectors $\{\bfx_{\alpha} \}_{\alpha \in A_+}$ in ${E_{\sbff \! \sbff}}^{-1} \{ \bfzero \}^{\perp}$ such that $E_{\sbff \! \sbff} \bfx_{\alpha} = \lambda_{\alpha} \bfx_{\alpha}$ where $\lambda_{\alpha} > 0$ for all $\alpha \in {\mathbb A_+}$.  There is also an orthonormal basis $\{\bfx_{\alpha}\}_{\alpha \in A_0}$ in ${E_{\sbff \! \sbff}}^{-1} \{ \bfzero \}$ with $E_{\sbff \! \sbff} \bfx_{\alpha} = \bfzero$ for all $\alpha \in A_0$.   This basis, which may be uncountable,  is automatically a basis of eigenvectors.   If we define $A = A_0 \cup A_+$ then $\{ \bfx_{\alpha}\}_{\alpha \in A}$ is a complete set of orthonormal eigenvectors in $H = {E_{\sbff \! \sbff}}^{-1} \{ \bfzero \} \oplus {E_{\sbff \! \sbff}}^{-1} \{ \bfzero \}^{\perp}$.  It follows that
\begin{eqnarray*}
\lefteqn{ \mbox{tr}(E_{\sbff \! \sbff}) = \mbox{$\sum_{\alpha \in A}$} \langle E_{\sbff \! \sbff}\bfx_{\alpha}, \bfx_{\alpha} \rangle =  \mbox{$\sum_{\alpha \in A_+}$} \langle E_{\sbff \! \sbff}\bfx_{\alpha}, \bfx_{\alpha} \rangle } \hspace{3cm} \\
& = & \mbox{$\sum_{\alpha \in A_+}$} \mbox{$\int_{\Omega}$}  | \langle \bfx_{\alpha}, \bff(\bomega) \rangle |^2 \mu(d \bomega) \\
& = & \mbox{$\int_{\Omega}$} \mbox{$\sum_{\alpha \in A_+}$} | \langle \bfx_{\alpha}, \bff(\bomega) \rangle |^2 \mu(d \bomega) \\
& = & \mbox{$\int_{\Omega}$} \|\bff(\bomega)\|^2 \mu(d \bomega) = {\mathbb E}[ \| \bff \|^2] = \dnorm \bff \dnorm^2 <  \infty. 
\end{eqnarray*}
Therefore $E_{\sbff \! \sbff}$ is nuclear and hence also compact \cite[p 279]{yos1}.  Note that
$$
\mbox{tr}(E_{\sbff \! \sbff}) = \mbox{$\sum_{\alpha \in A}$} \langle E_{\sbff \! \sbff} \bfx_{\alpha}, \bfx_{\alpha} \rangle = \mbox{$\sum_{\alpha \in A_+}$}\! \langle \lambda_{\alpha} \bfx_{\alpha}, \bfx_{\alpha} \rangle = \mbox{$\sum_{\alpha \in A_+}$}\! \lambda_{\alpha}.
$$
Consequently the operators $E_{\sbff \! \sbff} \in {\mathcal B}(H)$ and $E_{\sbfg \sbff} \in {\mathcal B}(H, K)$ satisfy the definitive properties
\begin{eqnarray}
\langle E_{\sbff \! \sbff}\, \bfx, \bfu \rangle & = & {\mathbb E}[ \langle \bfx, \bff \rangle \langle \bff, \bfu \rangle],  \label{covopp1} \\
\langle E_{\sbfg \sbff}\, \bfx, \bfy \rangle & = & {\mathbb E}[ \langle \bfx, \bff \rangle \langle \bfg, \bfy \rangle], \label{covopp2} \\
\mbox{tr}(E_{\sbff \! \sbff}) & = & {\mathbb E}[ \| \bff \|^2]  \label{covopp3}
\end{eqnarray}
for all $\bfx, \bfu \in H$ and $\bfy \in K$.  Thus we can regard these operators as covariance operators.  

\section{The generalized inverse auto-covariance operator}
\label{giac}

In this section we describe the generalized inverse auto-covariance operator.   We use an orthonormal basis of eigenvectors to construct a Fourier series representation of the auto-covariance $E_{\sbff \! \sbff}$ and hence define the generalized inverse auto-covariance ${E_{\sbff \! \sbff}}^{\dag}$.  We establish the important properties and pay particular attention to the general case where ${E_{\sbff \! \sbff}}^{\dag}$ is unbounded, closed, densely defined and self-adjoint.

Let $\{ \bfx_{\alpha} \}_{\alpha \in A}$ be a complete set of orthonormal eigenvectors for $E_{\sbff \! \sbff}$ in $H$ with corresponding eigenvalues $\{ \lambda_{\alpha} \}_{\alpha \in A}$.  The set $A_+ = \{ \alpha \mid \lambda_{\alpha} > 0 \}$ is at most a countable set but the set $A_0 = A \setminus A_+$ may be uncountable.  For each $\bfx \in H$ write
\begin{eqnarray*}
\bfx & = & \bfx_0 + \bfx_+ \\
& = & \mbox{$\sum_{\alpha \in A_0}$} \langle \bfx, \bfx_{\alpha} \rangle \bfx_{\alpha}  + \mbox{$\sum_{\alpha \in A_+}$} \langle \bfx, \bfx_{\alpha} \rangle \bfx_{\alpha} \in {E_{\sbff \! \sbff}}^{-1} \{ \bfzero\} \oplus {E_{\sbff \! \sbff}}^{-1}\{\bfzero\}^{\perp}
\end{eqnarray*}
and define a corresponding element $\bfu = E_{\sbff \! \sbff} \bfx \in E_{\sbff \! \sbff}(H)$ by the formula
$$
\bfu =  E_{\sbff \! \sbff} \bfx_0 + E_{\sbff \! \sbff} \bfx_+  =  E_{\sbff \! \sbff} \bfx_+  =  \mbox{$\sum_{\alpha \in A_+}$} \lambda_{\alpha} \langle \bfx, \bfx_{\alpha} \rangle \bfx_{\alpha}.
$$
Therefore $\bfu = \sum_{\alpha \in A_+} \langle \bfu, \bfx_{\alpha} \rangle \bfx_{\alpha}$ with $\langle \bfu, \bfx_{\alpha} \rangle = \lambda_{\alpha} \langle \bfx, \bfx_{\alpha} \rangle$ for each $\alpha \in A_+$ and so $\sum_{\alpha \in A_+} {\lambda_{\alpha}}^{-2} |\langle \bfu, \bfx_{\alpha} \rangle |^2 = \sum_{\alpha \in A_+} | \langle \bfx, \bfx_{\alpha} \rangle |^2 = \| \bfx_+ \|^2 < \infty$.
Conversely, if we are given $\bfu = \sum_{\alpha \in A_+} \langle \bfu, \bfx_{\alpha} \rangle \bfx_{\alpha}$ with $\sum_{\alpha \in A_+} {\lambda_{\alpha}}^{-2} | \langle \bfu, \bfx_{\alpha} \rangle |^2 < \infty$ then we can define $\bfx = \sum_{\alpha \in A_+} \lambda_{\alpha} \langle \bfu, \bfx_{\alpha} \rangle \bfx_{\alpha} \in H$ so that $E_{\sbff \! \sbff} \bfx = \bfu$.  Therefore $\bfu \in E_{\sbff \! \sbff}(H)$.  It follows that
$$
E_{\sbff \! \sbff}(H) = \{ \bfu \in H \mid \mbox{$\sum_{\alpha \in A_+}$}\, {\lambda_{\alpha}}^{-2} | \langle \bfu, \bfx_{\alpha} \rangle |^2 < \infty\} \subseteq {{E_{\sbff \! \sbff}}^{-1} \{ \bfzero \}}^{\perp}.
$$
There are two cases to consider.  If the index set $A_+$ is finite then for some $m \in {\mathbb N}$ we can write $A_+ = \{ j \in {\mathbb N} \mid j \leq m\}$.  In this case $E_{\sbff \! \sbff}(H) = {E_{\sbff \! \sbff}}^{-1}\{\bfzero\}^{\perp}$ is finite dimensional and closed, and the problem has already been solved \cite{fom1,how2}.  Henceforth we assume that $A_+$ is infinite and write $A_+ = {\mathbb N}$ with eigenvectors $\{\bfx_j\}_{j \in {\mathbb N}}$ and corresponding eigenvalues $\{ \lambda_j\}_{j \in {\mathbb N}}$ ordered in such a way that $\lambda_j \geq \lambda_{j+1} > 0$.  Now let $D({E_{\sbff \! \sbff}}^{\dag}) = E_{\sbff \! \sbff}(H) \oplus {E_{\sbff \! \sbff}}^{-1} \{ \bfzero \}$ and define ${E_{\sbff \! \sbff}}^{\dag}: D({E_{\sbff \! \sbff}}^{\dag}) \rightarrow {{E_{\sbff \! \sbff}}^{-1} \{ \bfzero \}}^{\perp}$ by setting
$$
{E_{\sbff \! \sbff}}^{\dag} \bfu = \mbox{$\sum_{j \in {\mathbb N}}$}\, {\lambda_j}^{-1} \langle \bfu, \bfx_j \rangle \bfx_j
$$
for each $\bfu \in D({E_{\sbff \! \sbff}}^{\dag})$.  We will use the above notation for the eigenvectors and eigenvalues throughout Section~\ref{giac} without further comment.

\subsection{The domain.}
\label{dgiac}

We will show that the domain $D({E_{\sbff \! \sbff}}^{\dag})$ is not closed.  Our definition of ${E_{\sbff \! \sbff}}^{\dag}$ is a natural definition.  If $\bfu \in E_{\sbff \! \sbff}(H)$ then there is a unique point $\bfx \in {E_{\sbff \!\sbff}}^{-1}(\{\bfzero\})^{\perp}$ such that $\bfu = E_{\sbff \! \sbff} \bfx$.  Hence we can define ${E_{\sbff \! \sbff}}^{\dag} \bfu = \bfx$.  If $\bfu \in {E_{\sbff \! \sbff}(H)}^{\perp}$ we define ${E_{\sbff \! \sbff}}^{\dag} \bfu = \bfzero$.  We begin by showing that $E_{\sbff \! \sbff}(H)$ is not closed.  We need to find $\{\bfu_n\}_{n \in {\mathbb N}} \subseteq E_{\sbff \! \sbff}(H)$ and $\bfu \notin E_{\sbff \! \sbff}(H)$ such that $\|\bfu_n - \bfu \| \rightarrow 0$ as $n \rightarrow \infty$.

To do this we need to construct a series $\sum_{j \in {\mathbb N}} \kappa_j$ that converges more slowly than $\sum_{j \in {\mathbb N}} \lambda_j$.  The following construction is taken from \cite{ash1}.  Define $\rho_j = \sum_{k=j}^{\infty} \lambda_k$ for each $j \in {\mathbb N}$ and define $\kappa_j = \lambda_j/\sqrt{\rho_j}$.  On the one hand $\kappa_j/\lambda_j = 1/ \sqrt{\rho_j} \rightarrow \infty$ as $j \rightarrow \infty$ and on the other hand
\begin{eqnarray*}
\mbox{$\sum_{j \in {\mathbb N}}$} \kappa_j & = & \mbox{$\sum_{j \in {\mathbb N}}$}\, (\rho_j - \rho_{j+1})/\sqrt{\rho_j} \\
& = & \mbox{$\sum_{j \in {\mathbb N}}$}\, (\sqrt{\rho_j} - \sqrt{\rho_{j+1}})(\sqrt{\rho_j} + \sqrt{\rho_{j+1}})/\sqrt{\rho_j} \\
& \leq & \mbox{$\sum_{j \in {\mathbb N}}$}\, 2(\sqrt{\rho_j} - \sqrt{\rho_{j+1}}) = 2\sqrt{\rho_1} < \infty.
\end{eqnarray*}
Thus $\sum_{j \in {\mathbb N}} \kappa_j$ is the desired series.  Since $\sum_{j\in {\mathbb N}} ({\kappa_j}/\kappa_1)^2 < \sum_{j\in {\mathbb N}} \kappa_j/\kappa_1 < \infty$ we can define $\bfu = \sum_{j \in {\mathbb N}} \kappa_j \bfx_j \in H$.  If we also define $\bfx_n = \sum_{j=1}^n (\kappa_j/\lambda_j) \bfx_j$ and $\bfu_n = E_{\sbff \! \sbff} \bfx_n = \sum_{j=1}^n \kappa_j \bfx_j \in E_{\sbff \! \sbff}(H)$ for each $n \in {\mathbb N}$ then $\| \bfu_n - \bfu \| \rightarrow 0$ as $n \rightarrow \infty$.  However $\{\bfx_n\}_{n\in {\mathbb N}}$ does not converge.  Therefore $\bfu \notin E_{\sbff \! \sbff}(H)$.  Equivalently we may say that $\bfu_n \in D({E_{\sbff \! \sbff}}^{\dag})$ with $\bfu_n \rightarrow \bfu \in H$ as $n \rightarrow \infty$ but with $\bfx_n = {E_{\sbff \! \sbff}}^{\dag} \bfu_n \in H$ for each $n \in {\mathbb N}$ such that $\{ \bfx_n \}_{n \in {\mathbb N}}$ diverges.  Thus $D({E_{\sbff \! \sbff}}^{\dag})$ is not closed.  

\subsection{The characteristic properties.}
\label{cpgiac}

We will show that ${E_{\sbff \! \sbff}}^{\dag}$ is unbounded, closed, densely defined and self-adjoint.

The operator ${E_{\sbff \! \sbff}}^{\dag}$ is unbounded because ${E_{\sbff \! \sbff}}^{\dag} \bfx_j = {\lambda_j}^{-1} \bfx_j$ for each $j \in {\mathbb N}$ with $\lambda_j \rightarrow 0$ as $j \rightarrow \infty$.

The following argument shows that ${E_{\sbff \! \sbff}}^{\dag}$ is closed.  Let $\{\bfu_n\}_{n \in {\mathbb N}} \subseteq D({E_{\sbff \! \sbff}}^{\dag})$.   Write $\bfu_n = \sum_{j \in {\mathbb N}} \langle \bfu_n, \bfx_j \rangle \bfx_j$ and ${E_{\sbff \! \sbff}}^{\dag} \bfu_n = \mbox{$\sum_{j \in {\mathbb N}}$}\, {\lambda_j}^{-1} \langle \bfu_n, \bfx_j \rangle \bfx_j$ for each $n \in {\mathbb N}$.  Now suppose that
$$
\| \bfu_n - \bfu \|^2 = \mbox{$\sum_{j \in {\mathbb N}}$} | \langle \bfu_n, \bfx_j \rangle - \langle \bfu, \bfx_j \rangle |^2  \rightarrow 0
$$
for some $\bfu \in H$ and that
$$
\| {E_{\sbff \! \sbff}}^{\dag} \bfu_n - \bfx \|^2 = \mbox{$\sum_{j \in {\mathbb N}}$} |  {\lambda_j}^{-1} \langle \bfu_n, \bfx_j \rangle - \langle \bfx, \bfx_j \rangle |^2 \rightarrow 0
$$
as $n \rightarrow \infty$ for some $\bfx \in H$.  Therefore
\begin{eqnarray*}
\mbox{$\sum_{j \in {\mathbb N}}$} | \langle \bfu_n, \bfx_j \rangle - \lambda_j \langle \bfx, \bfx_j \rangle |^2 & = & \mbox{$\sum_{j \in {\mathbb N}}$} {\lambda_j}^2 | {\lambda_j}^{-1} \langle \bfu_n, \bfx_j \rangle - \langle \bfx, \bfx_j \rangle |^2 \\
& \leq & {\lambda_1}^2 \mbox{$\sum_{j \in {\mathbb N}}$} | {\lambda_j}^{-1} \langle \bfu_n, \bfx_j \rangle - \langle \bfx, \bfx_j \rangle |^2 \rightarrow 0
\end{eqnarray*}
as $n \rightarrow \infty$.  Hence $\bfu_n \rightarrow \bfu = \sum_{j \in {\mathbb N}} \lambda_j \langle \bfx, \bfx_j \rangle \bfx_j$.  Now we have ${E_{\sbff \! \sbff}}^{\dag} \bfu = \bfx$ as required.  Thus ${E_{\sbff \! \sbff}}^{\dag}$ is a closed operator.

We show that $D({E_{\sbff \! \sbff}}^{\dag})$ is dense in $H$.  For each $\bfu = \sum_{\alpha \in A} \langle \bfu, \bfx_{\alpha} \rangle \bfx_{\alpha} \in H$ we can define a sequence $\{\bfu_n\}_{n \in {\mathbb N}} \subseteq D({E_{\sbff \! \sbff}}^{\dag})$ by setting
\begin{equation}
\label{un}
\bfu_n = \mbox{$\sum_{j=1}^n$} \langle \bfu, \bfx_j \rangle \bfx_j + \mbox{$\sum_{\alpha \in A_0}$} \langle \bfu, \bfx_{\alpha} \rangle \bfx_{\alpha} 
\end{equation}
such that $\| \bfu_n - \bfu \|^2 = \sum_{j=n+1}^{\infty} |\langle \bfu, \bfx_j \rangle |^2 \rightarrow 0$ as $n \rightarrow \infty$.  Thus ${E_{\sbff \! \sbff}}^{\dag}$ is densely defined. 

Finally we show that ${E_{\sbff \! \sbff}}^{\dag}$ is self-adjoint.  Suppose $\bfu, \bfv \in D({E_{\sbff \! \sbff}}^{\dag})$.  If we write $\bfu = \mbox{$\sum_{\alpha \in A}$} \langle \bfu, \bfx_{\alpha} \rangle \bfx_{\alpha}$ and $\bfv = \mbox{$\sum_{\alpha \in A}$} \langle \bfv, \bfx_{\alpha} \rangle \bfx_{\alpha}$ then we have
$$
{E_{\sbff \! \sbff}}^{\dag} \bfu = \mbox{$\sum_{j \in {\mathbb N}}$} {\lambda_j}^{-1} \langle \bfu, \bfx_j \rangle \bfx_j \quad  \mbox{and} \quad {E_{\sbff \! \sbff}}^{\dag} \bfv = \mbox{$\sum_{j \in {\mathbb N}}$} {\lambda_j}^{-1} \langle \bfv, \bfx_j \rangle \bfx_j.
$$
Consequently 
\begin{eqnarray*}
\langle {E_{\sbff \! \sbff}}^{\dag} \bfu, \bfv \rangle & = & \mbox{$\sum_{j \in {\mathbb N}}$} {\lambda_j}^{-1} \langle \bfu, \bfx_j \rangle \langle \bfx_j, \bfv \rangle \\
& = & \mbox{$\sum_{j \in {\mathbb N}}$} {\lambda_j}^{-1}\, \overline{\langle \bfv, \bfx_j \rangle \, \langle \bfx_j, \bfu \rangle} =  \overline{\langle {E_{\sbff \! \sbff}}^{\dag} \bfv, \bfu \rangle} =  \langle \bfu, {E_{\sbff \! \sbff}}^{\dag} \bfv \rangle.
\end{eqnarray*}
Thus ${E_{\sbff \! \sbff}}^{\dag}$ is self-adjoint.

\subsection{The standard properties.}
\label{spgiac} 

We justify our definitions by showing that ${E_{\sbff \! \sbff}}^{\dag}$ satisfies the standard properties associated with a generalized inverse operator.   Let $\bfu = \sum_{\alpha \in A} \langle \bfu, \bfx_{\alpha} \rangle \bfx_{\alpha} \in H$ and let $\{ \bfu_n \}_{n \in {\mathbb N}} \subseteq H$ be the sequence defined above in (\ref{un}) with $\bfu_n \in D({E_{\sbff \! \sbff}}^{\dag})$ for all $n \in {\mathbb N}$ and $\| \bfu_n - \bfu \| \rightarrow 0$ as $n \rightarrow \infty$.  Since $E_{\sbff \! \sbff} {E_{\sbff \! \sbff}}^{\dag}\bfu_n = \mbox{$\sum_{j \in {\mathbb N}}$} \langle \bfu_n, \bfx_j \rangle \bfx_j$ we can define 
$$
E_{\sbff \! \sbff}{E_{\sbff \! \sbff}}^{\dag}\bfu = \lim_{n \rightarrow \infty} \mbox{$\sum_{j \in {\mathbb N}}$} \langle \bfu_n, \bfx_j \rangle \bfx_j = \mbox{$\sum_{j \in {\mathbb N}}$} \langle \bfu, \bfx_j \rangle \bfx_j \in {E_{\sbff \! \sbff}}^{-1}\{ \bfzero\}^{\perp}.
$$
Therefore $\langle E_{\sbff \! \sbff}{E_{\sbff \! \sbff}}^{\dag}\bfu, \bfv \rangle = \mbox{$\sum_{j \in {\mathbb N}}$} \langle \bfu, \bfx_j \rangle \langle \bfx_j, \bfv \rangle$ for each $\bfu, \bfv \in H$.  For each $\bfv \in H$ we have $E_{\sbff \! \sbff} \bfv = \sum_{j \in {\mathbb N}} \lambda_j \langle \bfv, \bfx_j \rangle \bfx_j \in E_{\sbff \! \sbff}(H) \subseteq D({E_{\sbff \! \sbff}}^{\dag})$.  It follows that ${E_{\sbff \! \sbff}}^{\dag}E_{\sbff \! \sbff} \bfv = \mbox{$\sum_{j \in {\mathbb N}}$} \langle \bfv, \bfx_j \rangle \bfx_j \in {E_{\sbff \! \sbff}}^{-1}\{ \bfzero\}^{\perp}$ and hence that
$$
\langle {E_{\sbff \! \sbff}}^{\dag}E_{\sbff \! \sbff} \bfv, \bfu \rangle = \mbox{$\sum_{j \in {\mathbb N}}$} \langle \bfv, \bfx_j \rangle \langle \bfx_j, \bfu \rangle
$$
for each $\bfu, \bfv \in H$.  A similar argument to that used in the previous section now shows that $\langle [{E_{\sbff \! \sbff}}^{\dag}E_{\sbff \! \sbff}]^* \bfu, \bfv \rangle = \langle {E_{\sbff \! \sbff}}^{\dag}E_{\sbff \! \sbff} \bfu, \bfv \rangle$ for all $\bfu, \bfv \in H$.

We can now see that the operator ${E_{\sbff \! \sbff}}^{\dag}: D({E_{\sbff \! \sbff}}^{\dag}) \rightarrow H$ has the following properties.
\begin{enumerate}
\item $E_{\sbff \! \sbff}{E_{\sbff \! \sbff}}^{\dag} E_{\sbff \! \sbff} = E_{\sbff \! \sbff} \in {\mathcal B}(H)$.
\item ${E_{\sbff \! \sbff}}^{\dag} E_{\sbff \! \sbff} {E_{\sbff \! \sbff}}^{\dag} = {E_{\sbff \! \sbff}}^{\dag}:D({E_{\sbff \! \sbff}}^{\dag}) \rightarrow H$.
\item $[E_{\sbff \! \sbff}{E_{\sbff \! \sbff}}^{\dag}]^* = E_{\sbff \! \sbff}{E_{\sbff \! \sbff}}^{\dag} \in {\mathcal B}(H)$.
\item $[{E_{\sbff \! \sbff}}^{\dag}E_{\sbff \! \sbff}]^* = {E_{\sbff \! \sbff}}^{\dag}E_{\sbff \! \sbff} \in {\mathcal B}(H)$.
\end{enumerate}

\subsection{Some specific identities.}
\label{sigiac}

Let $T \in {\mathcal B}(H,K)$ and suppose that $T^* \bfy \in D({E_{\sbff \!\sbff}}^{\dag})$ and that $\bff(\bomega) \in D({E_{\sbff \!\sbff}}^{\dag})$ for $\mu$-almost all $\bomega \in \Omega$.  We have
\begin{eqnarray*}
E_{\sbff \! \sbff}({E_{\sbff \! \sbff}}^{\dag} T^* \bfy) & =  & \mbox{$\int_{\Omega}$} \langle {E_{\sbff \! \sbff}}^{\dag}T^* \bfy, \bff(\bomega) \rangle \bff(\bomega) \mu(d \bomega) \\
& =  & \mbox{$\int_{\Omega}$} \langle T^* \bfy, {E_{\sbff \! \sbff}}^{\dag} \bff(\bomega) \rangle  \bff(\bomega) \mu(d \bomega) \\
& =  & \mbox{$\int_{\Omega}$} \langle \bfy, T {E_{\sbff \! \sbff}}^{\dag} \bff(\bomega) \rangle  \bff(\bomega) \mu(d \bomega)
\end{eqnarray*}
because ${E_{\sbff \! \sbff}}^{\dag}$ is self-adjoint. Therefore
$$
T{E_{\sbff \! \sbff}}^{\dag} E_{\sbff \! \sbff} {E_{\sbff \! \sbff}}^{\dag}T^* \bfy = \mbox{$\int_{\Omega}$} \langle \bfy, T {E_{\sbff \! \sbff}}^{\dag} \bff(\bomega) \rangle T{E_{\sbff \! \sbff}}^{\dag} \bff(\bomega) \mu(d \bomega) = E_{\sbfk \sbfk} \bfy
$$
where we have written $\bfk = T {E_{\sbff \! \sbff}}^{\dag} \bff$ for convenience.   Therefore we have $E_{\sbfk \sbfk} = T{E_{\sbff \! \sbff}}^{\dag} E_{\sbff \! \sbff} {E_{\sbff \! \sbff}}^{\dag}T^* = T {E_{\sbff \! \sbff}}^{\dag} T^*$.  Similar arguments can be used to show that $E_{\sbfg \sbfk} = E_{\sbfg \sbff} {E_{\sbff \! \sbff}}^{\dag}T^*$ and $E_{\sbfk \sbfg} =  T {E_{\sbff \! \sbff}}^{\dag} E_{\sbff \sbfg}$.  The proof of the main result makes use of these specific identities.

\section{The null spaces of the covariance operators}
\label{nsco}

The next two results are important to the solution of Problem~\ref{p1}.   We show that the null space of $E_{\sbff \! \sbff}$ is a subspace of the null space of $E_{\sbfg \sbff}$ and hence deduce that $E_{\sbfg \sbff} = E_{\sbfg \sbff}{E_{\sbff \! \sbff}}^{\dag}E_{\sbff \! \sbff}$.

\begin{lem}
\label{nullspace}
Let $P = {E_{\sbff \! \sbff}}^{-1}\{\bfzero\}$ and $Q = {E_{\sbfg \sbff}}^{-1}\{\bfzero\}$ denote the null spaces of $E_{\sbff \! \sbff}$ and $E_{\sbfg \sbff}$ respectively.  Then $P \subseteq Q \subseteq H$.  $\hfill \Box$
\end{lem}

{\bf Proof} \quad  Let $\bfu \in P$. Then
$$
{\mathbb E}[ | \langle \bfu, \bff \rangle |^2 ] = {\mathbb E}[ \langle \bfu, \bff \rangle \langle \bff, \bfu \rangle ] = \langle E_{\sbff \! \sbff} \bfu, \bfu \rangle = 0.
$$
For each $\bfv \in K$ it follows that
$$
| \langle E_{\sbfg \sbff} \bfu, \bfv \rangle | =  | {\mathbb E}[ \langle \bfu, \bff \rangle \langle \bfg, \bfv \rangle ] | \leq {\mathbb E}[ | \langle \bfu, \bff \rangle |^2]^{1/2} {\mathbb E}[ | \langle \bfg, \bfv \rangle |^2] ^{1/2} = 0.
$$
Therefore $E_{gf} \bfu = \bfzero$.  Hence $\bfu \in Q$. $\hfill \Box$

\begin{cor}
\label{ginvprod}
Let $H, K$ be Hilbert spaces with $\bff \in L^2(\Omega,H)$ and $\bfg \in L^2(\Omega, K)$.  We have
$$
E_{\sbfg \sbff}(I - {E_{\sbff \! \sbff}}^{\dag} E_{\sbff \! \sbff}) = 0 \iff E_{\sbfg \sbff} = E_{\sbfg \sbff}{E_{\sbff \! \sbff}}^{\dag}E_{\sbff \! \sbff}.
$$
$\hfill \Box$
\end{cor}

{\bf Proof} \quad Let $\bfx \in H$ and write $\bfx = \sum_{j \in {\mathbb N}} \langle \bfx, \bfx_j \rangle \bfx_j + \sum_{\alpha \in A \setminus A_+} \langle \bfx, \bfx_{\alpha} \rangle \bfx_{\alpha}$.  We know that 
$$
{E_{\sbff \! \sbff}}^{\dag}E_{\sbff \! \sbff} \bfx = \mbox{$\sum_{j \in {\mathbb N}}$} \langle \bfx, \bfx_j \rangle \bfx_j.
$$
Therefore
$$
(I - {E_{\sbff \! \sbff}}^{\dag}E_{\sbff \! \sbff}) \bfx = \mbox{$\sum_{\alpha \in A \setminus A_+}$} \langle \bfx, \bfx_{\alpha} \rangle \bfx_{\alpha} \in {E_{\sbff \! \sbff}}^{-1}\{ \bfzero\}
$$
for all $\bfx \in H$ from which it follows that $E_{\sbfg \sbff} (I - {E_{\sbff \! \sbff}}^{\dag}E_{\sbff \! \sbff}) = 0$.  $\hfill \Box$
 
\section{Solution of the general estimation problem}
\label{sgep}

Let us return to the original problem.  Let $\bff \in L^2(\Omega, H)$ and $\bfg \in L^2(\Omega, K)$ be random functions with zero means.  We wish to find a closed, densely defined, linear operator $X:D(X) \subseteq H \rightarrow K$, a proximate observable function $\bff_{\epsilon}$ for each $\epsilon > 0$, with ${\mathbb E}[ \| \bff_{\epsilon} - \bff \|^2] < \epsilon$ and $\bff_{\epsilon}(\bomega) \in D(X)$ for $\mu$-almost all $\bomega \in \Omega$, and a corresponding estimate $\widehat{\bfg}_{\epsilon} = X\!\bff_{\epsilon}$ such that the mean square error ${\mathbb E}[ \| X\! \bff_{\epsilon} - \bfg \|^2]$ is minimized.

Suppose $\bff(\bomega) \in D({E_{\sbff\!\sbff}}^{\dag})$ for $\mu$-almost all $\bomega \in \Omega$ and let $X:D({E_{\sbff\!\sbff}}^{\dag}) \subseteq H \rightarrow K$ be defined by $X = T{E_{\sbff \! \sbff}}^{\dag}$ for some $T \in {\mathcal B}(H,K)$.  Take $\bff_{\epsilon} = \bff$ and let $\bfr = X\bff - \bfg$ = $T {E_{\sbff\!\sbff}}^{\dag} \bff - \bfg = \bfk - \bfg$.  Now
\begin{eqnarray*}
{\mathbb E}[ \| \bfr \|^2 ] & = & \mbox{tr}(E_{\sbfr \sbfr}) \\
& = & \mbox{tr}( E_{\sbfk \sbfk} - E_{\sbfg \sbfk} - E_{\sbfk \sbfg} + E_{\sbfg \sbfg}) \\
& = &  \mbox{tr} ( T {E_{\sbff \! \sbff}}^{\dag} T^* - E_{\sbfg \sbff} {E_{\sbff \! \sbff}}^{\dag} T^* - T {E_{\sbff \! \sbff}}^{\dag} E_{\sbff \sbfg} + E_{\sbfg \sbfg} ) \\
& = &  \mbox{tr} ( (T - E_{\sbfg \sbff} ){E_{\sbff \! \sbff}}^{\dag} E_{\sbff \! \sbff} {E_{\sbff \! \sbff}}^{\dag}(T^* -  E_{\sbff \sbfg})) + \mbox{tr} ( E_{\sbfg \sbfg} - E_{\sbfg \sbff} {E_{\sbff \! \sbff}}^{\dag} E_{\sbff \sbfg} ) \\
& = & \mbox{tr} ( E_{\sbfv \sbfv} )+ \mbox{tr} ( E_{\sbfg \sbfg} - E_{\sbfg \sbff} {E_{\sbff \! \sbff}}^{\dag} E_{\sbff \sbfg} ) 
\end{eqnarray*}
where we have written $\bfv = (T - E_{\sbfg \sbff}) {E_{\sbff \! \sbff}}^{\dag} \bff \in L^2(\Omega,K)$.  Therefore
\begin{eqnarray*}
{\mathbb E}[ \| \bfr \|^2 ] & = & {\mathbb E}[ \|\bfv\|^2] + \mbox{tr} ( E_{\sbfg \sbfg} - E_{\sbfg \sbff} {E_{\sbff \! \sbff}}^{\dag} E_{\sbff \sbfg} ) \\
& = & {\mathbb E}[ \| (T - E_{\sbfg \sbff}){E_{\sbff \! \sbff}}^{\dag} \bff \|^2] + \mbox{tr} ( E_{\sbfg \sbfg} - E_{\sbfg \sbff} {E_{\sbff \! \sbff}}^{\dag} E_{\sbff \sbfg} ).
\end{eqnarray*}
Thus the minimum occurs when $(T - E_{\sbfg \sbff}) {E_{\sbff \! \sbff}}^{\dag} \bff(\bomega) = \bfzero$ for $\mu$-almost all $\bomega \in \Omega$.  Hence we choose $T = E_{\sbfg \sbff} + B(I - E_{\sbff \! \sbff} {E_{\sbff \! \sbff}}^{\dag})$ where $B \in {\mathcal B}(H, K)$ is arbitrary.  Therefore $X = E_{\sbfg \sbff} {E_{\sbff \! \sbff}}^{\dag}$.   The minimum value of the expected mean-square error is
$$
{\mathbb E}[ \| E_{\sbfg \sbff} {E_{\sbff \! \sbff}}^{\dag} \bff - \bfg \|^2] = \mbox{tr} ( E_{\sbfg \sbfg} - E_{\sbfg \sbff} {E_{\sbff \! \sbff}}^{\dag} E_{\sbff \sbfg} ).
$$
Since $X = E_{\sbfg \sbff} {E_{\sbff \! \sbff}}^{\dag}$ we may assume $D(X) = D({E_{\sbff\!\sbff}}^{\dag})$.  Therefore $X$ is closed and densely defined.

Now suppose there is a set $S$ with $\mu(S) > 0$ and $\bff(\bomega) \notin D({E_{\sbff \! \sbff}}^{\dag})$ for $\bomega \in S$.  Let $\{ \bfx_j\}_{j \in {\mathbb N}}$ be a complete set of orthonormal eigenvectors for $E_{\sbff\!\sbff}$ in ${E_{\sbff\!\sbff}}^{-1}( \{ \bfzero\})^{\perp}$.  Let $n \in {\mathbb N}$ and define a proximate observable function $\bfp = \bff_n$ by setting $\bff_n(\bomega) = \sum_{j=1}^n \langle \bfx_j, \bff(\bomega) \rangle \bfx_j$ for each $\bomega \in \Omega$.  Thus
$$
\langle \bfx_j, \bfp(\bomega) \rangle = \left\{ \begin{array}{cl}
\langle \bfx_j, \bff(\bomega) \rangle & \mbox{for}\ j \leq n\\
0 & \mbox{otherwise}. \end{array} \right.
$$
Since $E_{\sbff\!\sbff} \bfx_j = \lambda_j \bfx_j$ it follows that
$$
E_{\sbff\!\sbff} (\mbox{$\sum_{j=1}^n$} {\lambda_j}^{-1} \langle \bfx_j, \bff(\bomega) \rangle \bfx_j ) = \mbox{$\sum_{j=1}^n$} \langle \bfx_j, \bff(\bomega \rangle \bfx_j = \bfp(\bomega)
$$
and so $\bfp(\bomega) \in E_{\sbff\!\sbff}(H) \subseteq D({E_{\sbff \! \sbff}}^{\dag})$ for all $\bomega \in \Omega$. Therefore the corresponding optimal estimate using $\bfp = \bff_n$ rather than $\bff$ is given by $\widehat{\bfg}_n = E_{\sbfg \sbfp}{E_{\sbfp \sbfp}}^{\dag}\bfp$ with error $\mbox{tr}(E_{\sbfg \sbfg} - E_{\sbfg \sbfp} {E_{\sbfp \sbfp}}^{\dag} E_{\sbfp \sbfg})$.  Now, for $j \leq n$, we have
$$
E_{\sbff\!\sbff} \bfx_j = \mbox{$\int_{\Omega}$} \langle \bfx_j, \bff(\bomega) \rangle \bfx_j \mu(d \bomega) = \mbox{$\int_{\Omega}$} \langle \bfx_j, \bfp(\bomega) \rangle \bfx_j \mu(d \bomega) = E_{\sbfp \sbfp} \bfx_j.
$$
Therefore $E_{\sbfp \sbfp} \bfx_j = \lambda_j \bfx_j$ for each $j \leq n$ and so ${E_{\sbfp \sbfp}}^{\dag} \bfx_j = {\lambda_j}^{-1} \bfx_j = {E_{\sbff \! \sbff}}^{\dag} \bfx_j$.  Now $E_{\sbfg \sbfp} \bfx_j = \int_{\Omega} \langle \bfx_j, \bfp(\bomega) \rangle \bfg(\bomega) \mu(d \bomega) = \int_{\Omega} \langle \bfx_j, \bff(\bomega) \rangle \bfg(\bomega) \mu(d \bomega) = E_{\sbfg \sbff} \bfx_j$ for $j \leq n$.  It follows, by linearity, that
\begin{eqnarray*}
E_{\sbfg \sbfp}{E_{\sbfp \sbfp}}^{\dag} \bfp(\bomega) & = & E_{\sbfg \sbfp} {E_{\sbff \! \sbff}}^{\dag} \bfp(\bomega) \\
& = & \mbox{$\sum_{j=1}^n$} \langle \bfx_j, \bff(\bomega) \rangle E_{\sbfg \sbfp} {E_{\sbff \! \sbff}}^{\dag} \bfx_j \\
& = & \mbox{$\sum_{j=1}^n$} {\lambda_j}^{-1} \langle \bfx_j, \bff(\bomega) \rangle E_{\sbfg \sbfp} \bfx_j \\
& = & \mbox{$\sum_{j=1}^n$} {\lambda_j}^{-1} \langle \bfx_j, \bff(\bomega) \rangle E_{\sbfg \sbff} \bfx_j \\
& = & E_{\sbfg \sbff} (\mbox{$\sum_{j=1}^n$} {\lambda_j}^{-1} \langle \bfx_j, \bff(\bomega) \rangle \bfx_j ) = E_{\sbfg \sbff} {E_{\sbff \! \sbff}}^{\dag} \bfp(\bomega)
\end{eqnarray*}
for all $\bomega \in \Omega$.  Since $\bfp = \bff_n$ the corresponding optimal estimate can now be written as $\widehat{\bfg}_n(\bomega) = E_{\sbfg \sbff} {E_{\sbff \! \sbff}}^{\dag} \bff_n(\bomega)$ for all $\bomega \in \Omega$.  Thus we may take $X = E_{\sbfg \sbff} {E_{\sbff \! \sbff}}^{\dag}$ as before.  The only difference is that we replace $\bff$ by $\bff_n$ for some suitably large value of $n \in {\mathbb N}$.  Note that ${\mathbb E}[ \|\bff_n - \bff\|^2 ] = \sum_{j=n+1}^{\infty} \int_{\Omega} | \langle \bfx_j, \bff(\bomega) \rangle |^2 \mu(d \bomega) \rightarrow 0$ as $n \rightarrow \infty$.

\section{A practical solution procedure}
\label{psp}

In practice we may be restricted to observing a projected component $\bfp(\bomega) = P \cdot \bff(\bomega)$ of the outcome $\bff(\bomega)$ where $P \in {\mathcal B}(H)$ is an orthogonal projection onto a closed subspace $M = P(H) \subseteq H$.  We would like to relate the restricted optimal estimate to the true optimal estimate.  

\begin{lem}
\label{gsrs}
Let $H, K$ be Hilbert spaces and let $P \in {\mathcal B}(H)$ be an orthogonal projection onto the closed subspace $M = P(H)$.  Let $\bff \in L^2(\Omega,H)$ and $\bfg \in L^2(\Omega,K)$ be zero-mean random functions with $\bfp(\bomega) = P \cdot \bff(\bomega) \in M$ and $\bfq(\bomega) = (I - P)\cdot \bff(\bomega) \in M^{\perp}$ the respective observable and unobservable components of $\bff(\bomega)$ for each $\bomega \in \Omega$.  If we define $\bfr = \bfq - E_{\sbfq \sbfp} {E_{\sbfp \sbfp}}^{\dag} \bfp$ we can rewrite the equation $X E_{\sbff\!\sbff} = E_{\sbfg \sbff}$ where $X =  E_{\sbfg \sbff}{E_{\sbff\!\sbff}}^{\dag}:D({E_{\sbff\!\sbff}}^{\dag}) \subseteq H \rightarrow K$ in the form
\begin{equation}
\label{lsealt}
\left[ \begin{array}{cc}
Y & Z \end{array} \right] \left[ \begin{array}{cc}
E_{\sbfp \sbfp} & 0 \\
E_{\sbfq \sbfp} & E_{\sbfr \sbfr} \end{array} \right] = \left[ \begin{array}{cc}
E_{\sbfg \sbfp} & E_{\sbfg \sbfr} \end{array} \right]
\end{equation}
where $Y:D({E_{\sbff\!\sbff}}^{\dag}) \cap M \subseteq M \rightarrow K$ and $Z:D({E_{\sbff\!\sbff}}^{\dag}) \cap M^{\perp} \subseteq M^{\perp} \rightarrow K$ are given by $Y = (E_{\sbfg \sbfp} - E_{\sbfg \sbfr}{E_{\sbfr \sbfr}}^{\dag}E_{\sbfq \sbfp}){E_{\sbfp \sbfp}}^{\dag}$ and $Z = E_{\sbfg \sbfr}{E_{\sbfr \sbfr}}^{\dag}$.   The optimal estimate for $\bfg$ is
\begin{equation}
\label{resest}
\widehat{\bfg} = E_{\sbfg \sbfp}{E_{\sbfp \sbfp}}^{\dag}\bfp + E_{\sbfg \sbfr}{E_{\sbfr\sbfr}}^{\dag} \bfr = \widehat{\bfg}_{\!M} + \bfr_{\!M}
\end{equation}
where $\widehat{\bfg}_{\!M} = E_{\sbfg \sbfg}{E_{\sbfp \! \sbfp}}^{\dag} \bfp$ is the restricted optimal estimate.  The components $\widehat{\bfg}_M$ and $\bfr_M$ are uncorrelated and the error in the restricted estimate is
\begin{equation}
\label{reserror}
{\mathbb E}[ \| \widehat{\bfg}_M - \bfg \|^2] = {\mathbb E}[ \| \widehat{\bfg} - \bfg \|^2] + \mbox{\rm tr}(E_{\sbfg \sbfr}{E_{\sbfr \sbfr}}^{\dag} E_{\sbfr \sbfg}).
\end{equation}
$\hfill \Box$
\end{lem}

{\bf Proof}\quad The equation $XE_{\sbff\!\sbff} = E_{\sbfg \sbff}$ is equivalent to the equation
$$
\left[ \begin{array}{cc}
\!\!Y\! &\! Z\!\! \end{array} \right] \cdot \left[ \begin{array}{cc}
\!\!E_{\sbfp \sbfp} & E_{\sbfp \sbfq}\!\! \\
\!\!E_{\sbfq \sbfp} & E_{\sbfq \sbfq}\!\! \end{array} \right] \left[ \begin{array}{cc}
\!I & \!\!- {E_{\sbfp\sbfp}}^{\dag} E_{\sbfp \sbfq} \!\!\\
\!0 & I \end{array} \right] = \left[ \begin{array}{cc}
\!\!E_{\sbfg \sbfp} \!&\! E_{\sbfg \sbfq}\!\! \end{array} \right] \left[ \begin{array}{cc}
\!I & \!\!- {E_{\sbfp\sbfp}}^{\dag} E_{\sbfp \sbfq} \!\!\\
\!0 & I \end{array} \right].
$$
If we evaluate the matrix products and use the identities $E_{\sbfq \sbfp} = E_{\sbfq \sbfp} {E_{\sbfp \sbfp}}^{\dag} E_{\sbfp \sbfp}$ and $E_{\sbfr \sbfr} =  E_{\sbfq \sbfq} - E_{\sbfq \sbfp}{E_{\sbfp\sbfp}}^{\dag}E_{\sbfp \sbfq}$ we obtain (\ref{lsealt}).   Solving $Z E_{\sbfr \sbfr} = E_{\sbfg \sbfr}$ gives $Z = E_{\sbfg \sbfr}{E_{\sbfr \sbfr}}^{\dag}$ and solving 
$Y E_{\sbfp\sbfp}E_{\sbfp \sbfq} + Z E_{\sbfq \sbfp} = E_{\sbfg \sbfp}$ gives $Y = E_{\sbfg \sbfp}{E_{\sbfp \sbfp}}^{\dag} - Z E_{\sbfq \sbfp} {E_{\sbfp \sbfp}}^{\dag}$.  Substituting for $Z$ shows that $Y = E_{\sbfg \sbfp}{E_{\sbfp \sbfp}}^{\dag} - E_{\sbfg \sbfr}{E_{\sbfr \sbfr}}^{\dag}  E_{\sbfq \sbfp} {E_{\sbfp \sbfp}}^{\dag}$ as required.  Hence
\begin{eqnarray*}
\widehat{\bfg} = Y \bfp + Z \bfq & = & E_{\sbfg \sbfp} {E_{\sbfp \sbfp}}^{\dag} \bfp + E_{\sbfg \sbfr}{E_{\sbfr \sbfr}}^{\dag}( \bfq - E_{\sbfq \sbfp}{E_{\sbfp \sbfp}}^{\dag} \bfp) \\
& = & E_{\sbfg \sbfp}{E_{\sbfp \sbfp}}^{\dag} \bfp + E_{\sbfg \sbfr}{E_{\sbfr\sbfr}}^{\dag} \bfr
\end{eqnarray*}
which is (\ref{resest}).  We note that $E_{\sbfp \sbfr} = E_{\sbfp \sbfq} - E_{\sbfp \sbfp}{E_{\sbfp \sbfp}}^{\dag} E_{\sbfp \sbfq} = 0$ which shows that the components $\widehat{\bfg}_{\!M}$ and $\bfr_{\!M}$ are uncorrelated.  We know from the previous section that $\widehat{\bfg} = E_{\sbfg \sbff}{E_{\sbff \! \sbff}}^{\dag} \bff$ and so (\ref{resest}) gives 
$$
E_{\sbfg \sbff} {E_{\sbff \! \sbff}}^{\dag} \bff = E_{\sbfg \sbfp}{E_{\sbfp \sbfp}}^{\dag} \bfp + E_{\sbfg \sbfr}{E_{\sbfr\sbfr}}^{\dag} \bfr.
$$
Therefore
$$
E_{\sbfg \sbff} {E_{\sbff \! \sbff}}^{\dag} E_{\sbff \sbfg} = E_{\sbfg \sbfp}  {E_{\sbfp \sbfp}}^{\dag} E_{\sbfp \sbfg} + E_{\sbfg \sbfr} {E_{\sbfr \sbfr}}^{\dag} E_{\sbfr \sbfg}.
$$
Now we can use this relationship and the known error estimates
$$
{\mathbb E}[ \| \widehat{\bfg} - \bfg \|^2] = \mbox{tr}( E_{\sbfg \sbfg} - E_{\sbfg \sbff} {E_{\sbff \! \sbff}}^{\dag} E_{\sbff \sbfg})
$$
and
$$
{\mathbb E}[ \| \widehat{\bfg}_M - \bfg \|^2] = \mbox{tr}( E_{\sbfg \sbfg} - E_{\sbfg \sbfp} {E_{\sbfp \sbfp}}^{\dag} E_{\sbfp \sbfg})
$$
to deduce (\ref{reserror}). $\hfill \Box$

\section{A hypothetical example}
\label{he}
\begin{figure}[htb]
\begin{center}
\includegraphics[width=8cm]{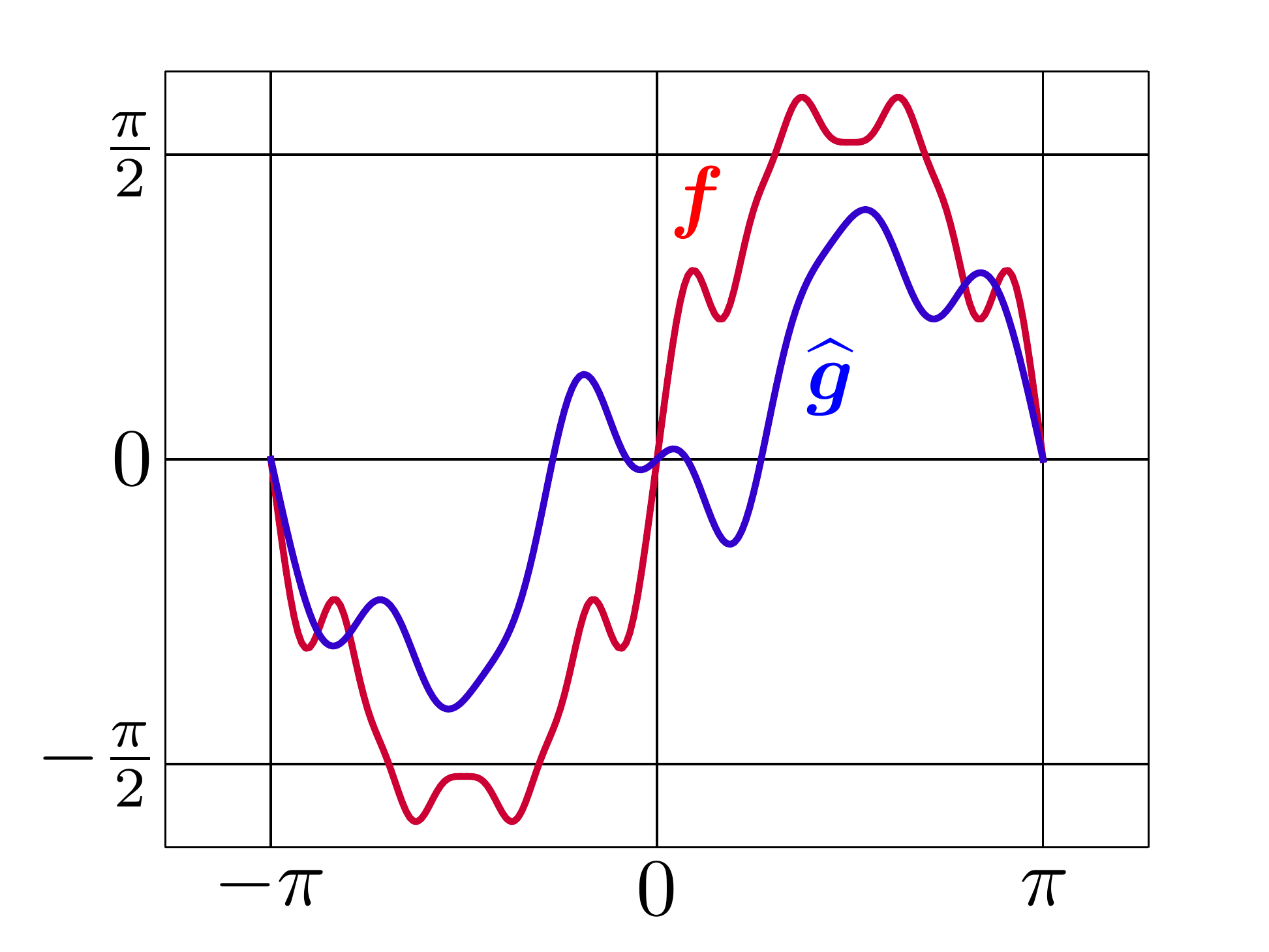}
\includegraphics[width=8cm]{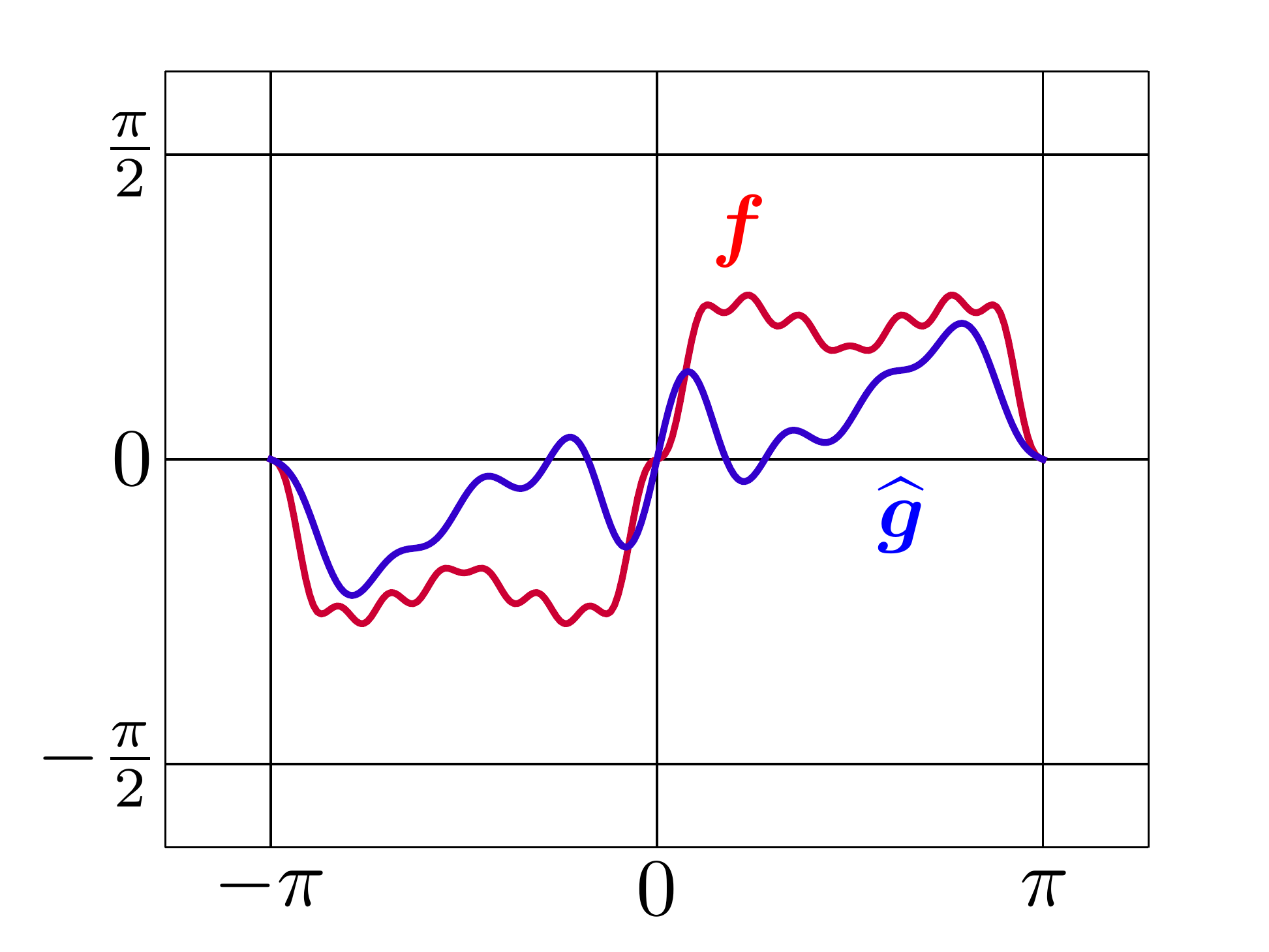}
\includegraphics[width=8cm]{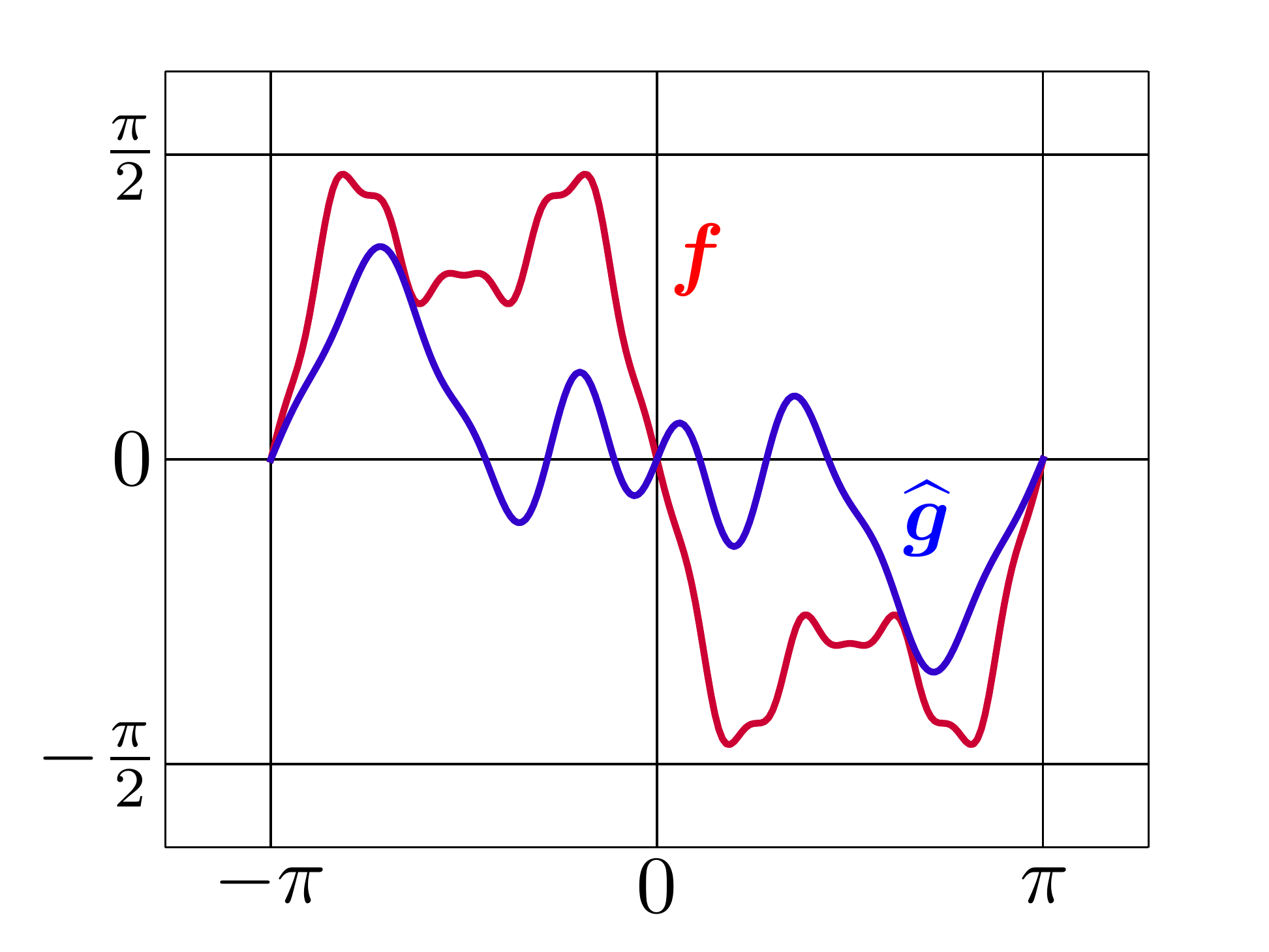}
\includegraphics[width=8cm]{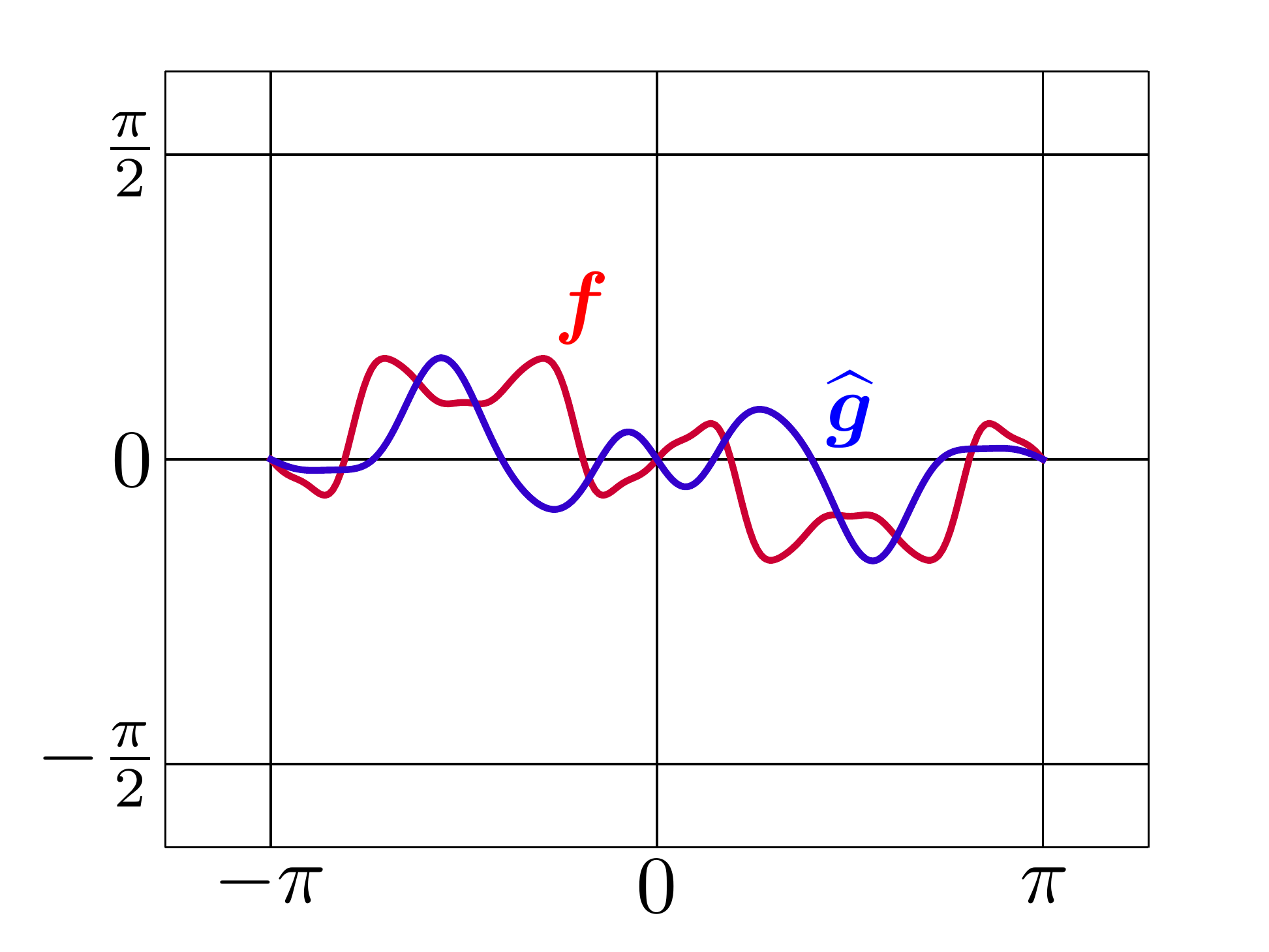}
\end{center}
\caption{The random function pairs $(\bff, \widehat{\bfg}) = (\bff, X\! \bff)$ for trials $2$ (top left), $4$ (top right), $5$ (bottom left) and $6$ (bottom right) showing a typical range of outcomes.  There is no estimation error in this example and so $\widehat{\bfg} = \bfg$ in each of these trials.}
\label{fig1}
\end{figure}

The functions $\varphi, \psi: (-\pi,\pi) \rightarrow {\mathbb R}$ defined by $\varphi(t) = \pi\, \mbox{sgn}(t)/4$ and $\psi(t) = t/2$ can be represented by the Fourier series
$$
\varphi(t) \sim  \mbox{$\sum_{k \in 2{\mathbb N}-1}$}\, \frac{1}{k} \sin kt \quad \mbox{and} \quad \psi(t) \sim \mbox{$\sum_{j \in {\mathbb N}}$} \frac{(-1)^{j+1}}{j} \sin jt. 
$$
Equivalently we may represent these functions as elements of the Hilbert space $\ell^{\,2}$ by the vectors
$$
\bphi \sim \left[ \begin{array}{c}
1 \\
0 \\
1/3 \\
0 \\
1/5 \\
\vdots \end{array} \right] \quad \mbox{and} \quad \bpsi \sim \left[ \begin{array}{c}
1 \\
-1/2 \\
1/3 \\
-1/4 \\
1/5 \\
\vdots \end{array} \right].
$$
Define a hypothetical experiment with outcomes $\bomega = \{ \omega_j \}_{j \in {\mathbb N}} \in \ell\! \stackrel{\infty}{\rule{0cm}{0.08cm}}$ where the coordinates $\omega_j \in {\mathbb R}$ for each $j \in {\mathbb N}$ are independent identically distributed random variables with cumulative distribution function $F:[-1,1] \rightarrow [0,1]$ defined by $F(t) = t/2 + 1/2$.  Let $\bff, \bfg \in L^2(\ell\! \stackrel{\infty}{\rule{0cm}{0.08cm}}, \ell^{\,2})$ be random functions with
$$
\bff = \left[ \begin{array}{c}
\omega_1 + \omega_2 \\
0 \\
(\omega_2 + \omega_3)/3 \\
0 \\
(\omega_3 + \omega_4)/5 \\
\vdots \end{array} \right] \quad \mbox{and} \quad \bfg = \left[ \begin{array}{c}
\omega_1 \\
- \omega_2/2 \\
\omega_3/3 \\
-\omega_4/4 \\
\omega_5/5 \\
\vdots \end{array} \right].
$$ 
The self-adjoint operator $E_{\sbff \! \sbff}$ can be represented by an infinite matrix $E_{\sbff \! \sbff} = [f\!\!f_{ij}]$ where $f\!\!f_{11} = 2/3$, $f\!\!f_{13} = 1/9$,
$$
f\!\!f_{k,k-2} = 1/[3k(k-2)], \ f\!\!f_{k,k} = 2/[3k^2], \ \mbox{and}\ f\!\!f_{k,k+2} = 1/[3k(k+2)]
$$
for each $k \in 2{\mathbb N}+1$, and $f\!\!f_{ij} = 0$ otherwise.  The operator $E_{\sbfg \sbff}$ can be represented by an infinite matrix $E_{\sbfg \sbff} = [g\!f_{ij}]$ where $g\!f_{11} = 1/3$,
$$
g\!f_{j, 2j-3} = (-1)^{j-1}\!/[3j(2j-3)]\ \mbox{and}\ g\!f_{j, 2j-1} = (-1)^{j-1}\!/[3j(2j-1)]
$$
for all $j \in {\mathbb N}+1$, and $g\!f_{ij} = 0$ otherwise.  Despite the structural simplicity of $E_{\sbff \! \sbff}$ it is a non-trivial task to calculate ${E_{\sbff \! \sbff}}^{\dag}$.  We can use elementary row operations to reduce the operator matrix to upper triangular form 
$$
U_{\sbff \! \sbff} = \left[ \begin{array}{cccccccccc}
2/3 & 0 & 1/9 & 0 & 0 & 0 & 0 & 0 & 0 & \cdots \\
0 & 0 & 0 & 0 & 0 & 0 & 0 & 0 & 0 & \cdots \\
0 & 0 & 1/18 & 0 & 1/45 & 0 & 0 & 0 & 0 & \cdots \\
0 & 0 & 0 & 0 & 0 & 0 & 0 & 0 & 0 & \cdots \\
0 & 0 & 0 & 0 & 4/225 & 0 & 1/105 & 0 & 0 & \cdots \\
0 & 0 & 0 & 0 & 0 & 0 & 0 & 0 & 0 & \cdots \\
0 & 0 & 0 & 0 & 0  & 0 & 5/588 & 0 & 1/189 & \cdots \\
0 & 0 & 0 & 0 & 0 & 0 & 0 & 0 & 0 & \cdots \\
0 & 0 & 0 & 0 & 0 & 0 & 0 & 0 & 2/405 & \cdots \\
\vdots & \vdots & \vdots & \vdots & \vdots & \vdots & \vdots & \vdots & \vdots & \ddots \end{array} \right]
$$ 
but a general formula for the elements on the leading diagonal is far from obvious.  We can gain some insight into the general calculation if we write
$$
E_{\sbff \! \sbff} = \bfc_1 {\bfc_1}^* + \bfc_3 {\bfc_3}^* + \bfc_5 {\bfc_5}^* + \cdots
$$
where we define $\bfc_{2j-1} = c_{2j-1,2j-1}\bfe_{2j-1} + c_{2j+1,2j-1}\bfe_{2j+1}$ for each $j \in {\mathbb N}$ and $\{\bfe_k\}_{k \in {\mathbb N}}$ are the standard basis vectors.   If we now equate coefficients we can see that ${c_{11}}^2 = f\!\!f_{11}$, $c_{31}c_{11} = f\!\!f_{31}$, ${c_{31}}^2 + {c_{33}}^2 = f\!\!f_{33}$, $c_{53}c_{33} = f\!\!f_{53}$, ${c_{53}}^2 + {c_{55}}^2 = f\!\!f_{55}$, $c_{75}c_{55} = f\!\!f_{75}$, ${c_{75}}^2 + {c_{77}}^2 = f\!\!f_{77}, \ldots$ and so on.  Solving these equations gives ${c_{11}}^2 = 2/3$, ${c_{31}}^2 = 1/54$, ${c_{33}}^2 = 1/18$, ${c_{53}}^2 = 2/225$, ${c_{55}}^2 = 4/225$, ${c_{75}}^2 = 1/196$, ${c_{77}}^2 = 5/588,\ldots$ and so on.  This suggests that the process actually defines the diagonal elements of the reduced matrix.  It turns out that it also defines the elementary row operations.  The coefficients $c_{2k+1, 2k+1}$ and $c_{2k+1,2k-1}$ are defined by the recursions
$$
{c_{2k+1,2k+1}}^2 = f\!\!f_{2k+1,2k+1} - \frac{{\\f\!\!f_{2k+1,2k-1}}^2}{{c_{2k-1,2k-1}}^2}\quad \mbox{and}\quad c_{2k+1,2k-1} = \frac{f\!\!f_{2k+1,2k-1}}{c_{2k-1,2k-1}}
$$
for each $k \in {\mathbb N}+1$ with $c_{11} = 2/3$ and $c_{31} = 1/9$.  If we define a sequence of lower triangular elementary operator matrices $L_k = [l_{k,ij}]$ by setting $l_{k,ii} = 1$, $l_{k,2k+1,2k-1} = (-1) c_{2k+1,2k-1}$ for each $k \in {\mathbb N}$ and $l_{k,ij} = 0$ otherwise, then we have
$$
L_{k-1} \cdots L_2 L_1 \cdot E_{\sbff \! \sbff} \cdot {L_1}^* {L_2}^* \cdots {L_{k-1}}^* = \left[ \begin{array}{cc}
D_{\sbff\! \sbff, [1,2k]} & 0 \\
 0 & E_{\sbff\! \sbff,\, [2k+1, \infty)} \end{array} \right] 
 $$
where $D_{\sbff \! \sbff,[1,2k]} = [d_{ij}] \in {\mathbb C}^{2k\times 2}$ is a diagonal matrix with
$$
d_{2\ell-1,2\ell-1} = {c_{2\ell-1,2\ell-1}}^2
$$
for $\ell \in {\mathbb N}$ and $d_{ij} = 0$ otherwise, and where $E_{\sbff\!\sbff, [2k+1, \infty)}$ denotes the operator matrix formed by deleting the first $2k$ rows and columns from $E_{\sbff\!\sbff}$. If we define $M_k = {L_k}^{-1}$ then it can be seen that
\begin{eqnarray*}
{E_{\sbff\!\sbff}}^{\dag} & = & M_1 M_2 \cdots M_{k-1}\left[ \begin{array}{cc}
{D_{\sbff\! \sbff, [1,2k]}}^{\dag} & 0 \\
 0 & {E_{\sbff\! \sbff,\, [2k+1, \infty)}}^{\dag} \end{array} \right] {M_{k-1}}^* \cdots {M_2}^*{M_1}^* \\
& = & \left[ \begin{array}{cc}
{E_{\sbff\! \sbff, [1,2k]}}^{\dag} & 0 \\
 0 & {E_{\sbff\! \sbff,\, [2k+1, \infty)}}^{\dag} \end{array} \right]
\end{eqnarray*}
for each $k \in {\mathbb N}$.  We know from the operator matrix representation of $E_{\sbff \! \sbff}$ that the trace is given by
$$
\mbox{tr}(E_{\sbff \! \sbff}) = 2/3 \left[ 1^2 + 1/3^2 + 1/5^2 + \cdots \right] = \pi^2/12 < \infty.
$$
Therefore $E_{\sbff \! \sbff}$ is a nuclear operator and hence ${E_{\sbff \! \sbff}}^{\dag}$ is closed and unbounded.  Some elementary algebra using M{\sc atlab} now suggests that we can represent the generalized inverse operator ${E_{\sbff\!\sbff}}^{\dag}$ in infinite matrix form as
$$
{E_{\sbff\!\sbff}}^{\dag} = \left[ \begin{array}{cccccccc}
3  & 0 & -9 & 0 & 15 & 0 & -21 & \cdots \\
0 & 0 & 0 & 0 & 0 & 0 & 0 & \cdots \\
-9 & 0 & 54 & 0 & -90 & 0 & 126 & \cdots \\
0 & 0 & 0 & 0 & 0 & 0 & 0 & \cdots \\
15 & 0 & -90 & 0 & 225 & 0 & -315 & \cdots \\
0 & 0 & 0 & 0 & 0 & 0 & 0 & \cdots \\
-21 & 0 & 126 & 0 & -315 & 0 & 588 & \cdots \\
\vdots & \vdots & \vdots & \vdots & \vdots & \vdots & \vdots & \ddots \end{array} \right] 
$$
where
$$
{f\!\!f_{2k-1, 2\ell-1}}^{\dag} = (-1)^{2k+2\ell-2}\cdot \min\{k,\ell\} \cdot 3(2k-1)(2\ell-1)
$$
for $k,\ell \in {\mathbb N}$ with ${f\!\!f_{ij}}^{\dag} = 0$ otherwise.   Now the matrix representation for $X = E_{\sbfg \sbff}{E_{\sbff\!\sbff}}^{\dag}$ is given by
$$
X = [x_{ij}] = \left[ \begin{array}{cccccccc}
1 & 0 & -3 & 0 & 5 & 0 & -7 & \cdots \vspace{0.1cm} \\
0 & 0 & \frac{-3}{2} & 0 & \frac{5}{2} & 0 & \frac{-7}{2} & \cdots \vspace{0.1cm} \\
0 & 0 & 0 & 0 & \frac{5}{3} & 0 & \frac{-7}{3} & \cdots \vspace{0.1cm} \\
0 & 0 & 0 & 0 & 0 & 0 & \frac{-7}{4} & \cdots \vspace{0.1cm} \\
0 & 0 & 0 & 0 & 0 & 0 & 0 & \cdots \\
\vdots & \vdots & \vdots & \vdots & \vdots & \vdots & \vdots & \ddots \end{array} \right]
$$
where $x_{11} = 1$ with
$$
x_{i,2\ell-1} = (-1)^{\ell-1} \cdot (2\ell -1)/i
$$
for $i < 2\ell-1$ and $i,\ell \in {\mathbb N}+1$ and $x_{ij} = 0$ otherwise.  The matrix representation for $X:D({E_{\sbff\!\sbff}}^{\dag}) \rightarrow K$ shows that it is unbounded and so we must be careful when calculating images for elements that are not in $D({E_{\sbff\!\sbff}}^{\dag})$.  Consider the calculation $\widehat{\bfg} = X\! \bff$.  Define
$$
\bff_{2n-1} = \sum_{k=1}^{n-1} (\omega_k + \omega_{k+1}) \bfe_{2k-1}/(2k-1) + \omega_n \bfe_{2n-1}/(2n-1)
$$
for each $n \in {\mathbb N}$.  Thus
$$
{\mathbb E}[\| \bff_{2n-1} - \bff \|^2] = (1/3)/(2n-1)^2 + (2/3) \sum_{k=n+1}^{\infty} 1/(2k-1)^2 \rightarrow 0
$$
as $n \rightarrow \infty$.  Now $\widehat{\bfg}_{2n-1} = X\! \bff_{2n-1}$ is given by 
$$
\widehat{\bfg}_{2n-1} = \left[ \begin{array}{c}
\sum_{k=1}^{n-1} (-1)^{k-1} (\omega_k + \omega_{k+1}) + (-1)^{n-1} \omega_n \vspace{0.1cm}\\
\sum_{k=2}^{n-1} (-1)^{k-1}\frac{(\omega_k + \omega_{k+1})}{2} + (-1)^{n-1} \frac{\omega_n}{2} \vspace{0.1cm}\\
\sum_{k=3}^{n-1} (-1)^{k-1}\frac{(\omega_k + \omega_{k+1})}{3} + (-1)^{n-1} \frac{\omega_n}{3} \\
\vdots \\
(-1)^{n-2}\frac{(\omega_{n-1} + \omega_n)}{(n-1)} + (-1)^{n-1} \frac{\omega_n}{(n-1)} \vspace{0.1cm}\\
(-1)^{n-1} \frac{\omega_n}{n} \\
0 \\
\vdots \end{array} \right] = \left[ \begin{array}{c}
\omega_1 \vspace{0.1cm}\\
-\frac{\omega_2}{2} \vspace{0.1cm}\\ 
\frac{\omega_3}{3} \vspace{0.1cm}\\
\vdots \\
(-1)^{n-2} \frac{\omega_{n-1}}{n-1} \vspace{0.1cm} \\
(-1)^{n-1} \frac{\omega_n}{n} \\
0 \\
\vdots
\end{array} \right]
$$
for all $m \in {\mathbb N}$.  This shows that
$$
{\mathbb E}[ \| \widehat{\bfg}_{2n-1} - \bfg \|^2] = (2/3) (1/(n+1)^2 + 1/(n+2)^2 + \cdots ) \rightarrow 0
$$
as $n \rightarrow \infty$ and so $\widehat{\bfg} = \bfg$.  For $k=8$ we have $E_{\sbfg \sbff} \approx E_{\sbfg \sbff,[1,8] \times [1,16]} \in {\mathbb C}^{8 \times 16}$ and ${E_{\sbff\!\sbff}}^{\dag} \approx {E_{\sbff\!\sbff,[1,16]}}^{\dag} \in {\mathbb C}^{16 \times 16}$ which gives
$$
X \approx E_{\sbfg \sbff,[1,8] \times [1,16]}{E_{\sbff\!\sbff,[1,16]}}^{\dag} = \left[ \begin{array}{cccccccc}
1 & 0 & -3 & 0 & 5 & \cdots & -15 & 0 \\
0 & 0 & \frac{-3}{2} & 0 & \frac{5}{2} & \cdots & \frac{-15}{2} & 0 \vspace{0.1cm} \\
0 & 0 & 0 & 0 & \frac{5}{3} & \cdots & \frac{-15}{3} & 0 \\
\vdots & \vdots & \vdots & \vdots & \vdots & \ddots & \vdots & \vdots \\
0 & 0 & 0 & 0 & 0 & \cdots & \frac{-15}{8} & 0 \end{array} \right]. 
$$
We performed ten trials.  The random function pairs $(\bff, \widehat{\bfg}) = (\bff, X\! \bff)$ for trials $2$, $4$, $5$ and $6$ are shown in F{\sc igure} \ref{fig1}. 

The trials used uniformly distributed pseudo-random numbers on $[-1,1]$ generated in M{\sc atlab}.  The results of trials $2$, $4$, $5$ and $6$ show a typical range of outcomes.  The corresponding pseudo-random numbers were
\begin{eqnarray*}
\bomega_2 & = & [0.9150, 0.9298, -0.6848, 0.9412, 0.9143, -0.0292, 0.6006, -0.7162] \\
\bomega_4 & = & [0.3575, 0.5155, 0.4863, -0.2155, 0.3110, -0.6576, 0.4121, -0.9363] \\
\bomega_5 & = & [-0.4462, -0.9077, -0.8057, 0.6469, 0.3897, -0.3658, 0.9004, -0.9311] \\
\bomega_6 & = & [ -0.1225, -0.2369, 0.5310, 0.5904, -0.6263, -0.0205, -0.1088, 0.2926].
\end{eqnarray*}
In this example it is easy to check that
$$
\mbox{tr}\left(E_{\sbfg\sbfg, [1,8]} - E_{\sbfg \sbff, [1,8] \times [1,16]} {E_{\sbff\!\sbff, [1, 16]}}^{\dag} E_{\sbff \sbfg, [1,16]\times [1,8]} \right) = 0
$$
and hence there is no estimation error.  We can explain this by noting that $\bff$ contains complete information about the outcome $\bomega$ and that we have used known theoretical information to construct the key matrices $E_{\sbff\!\sbff}$ and $E_{\sbfg \sbff}$.  In addition there are no observation errors in our model.  In practice $\bff$ may not contain complete information about the outcome, the observed values of $\bff(\bomega)$ will normally contain measurement errors, and the key matrices will likely be estimated from experimental data obtained under laboratory conditions where both $\bff(\bomega)$ and $\bfg(\bomega)$ can be observed. 

\section{Conclusions and future research}
\label{conc}

We have shown that the optimal least squares linear filter can be extended to estimation of random functions with values in infinite-dimensional Hilbert spaces.  In particular we have shown that in those instances where the generalized inverse auto-covariance is an unbounded linear operator it is nevertheless closed and densely defined.  Our future research will consider applications to signal processing and possible applications to the inversion of linear operator pencils where the resolvent operator has an isolated essential singularity at the origin \cite{alb2}.  These operators may arise in input retrieval problems for infinite-dimensional linear control systems \cite[Section 8.4.1, pp 261\textendash 262]{avr1} or in the solution of infinite systems of ordinary differential equations \cite[Section 8]{alb1}.

\end{document}